\documentclass[11pt,reqno]{amsart}
\usepackage[T1]{fontenc}
\usepackage{lmodern}
\usepackage[a4paper,margin=27mm]{geometry}
\usepackage{amsmath,amssymb,mathtools}
\usepackage{graphicx}
\usepackage{microtype}
\usepackage{enumitem}
\usepackage{booktabs}
\usepackage{array}
\usepackage{longtable}
\usepackage{tikz}
\usetikzlibrary{arrows.meta,positioning,calc}
\usepackage{pgfplots}
\pgfplotsset{compat=1.18}
\usepgfplotslibrary{fillbetween}
\usepackage[backend=biber,style=numeric,sorting=nyt,giveninits=true,isbn=false,eprint=false,url=false]{biblatex}
\DeclareFieldFormat{postnote}{#1}
\DeclareFieldFormat{multipostnote}{#1}
\DeclareSourcemap{
    \maps[datatype=bibtex]{
        \map{
            \step[fieldset=URL, null]
            \step[fieldset=DOI, null]
            \step[fieldset=MRREVIERER, null]
            \step[fieldset=ISBN, null]
            \step[fieldset=ISSN, null]
            \step[fieldset=SERIES, null]
        }
    }
}
\usepackage[hypertexnames=false,colorlinks=true,linkcolor=blue!55!black,citecolor=blue!55!black,urlcolor=blue!55!black]{hyperref}
\hypersetup{pdftitle={The sharp Hölder constant in Mori's theorem (restructured exposition)},pdfauthor={Qiang Gao}}
\numberwithin{equation}{section}
\newtheorem{theorem}{Theorem}[section]
\newtheorem{lemma}[theorem]{Lemma}
\newtheorem{proposition}[theorem]{Proposition}
\newtheorem{corollary}[theorem]{Corollary}
\theoremstyle{definition}

\theoremstyle{remark}
\newtheorem{remark}[theorem]{Remark}

\newcommand{\set}[1]{\left \{ #1 \right \}}
\makeatletter
\newcommand{\norm@auto}[2]{%
    \begingroup
    \setbox0=\hbox{$\m@th#1#2$}%
    \setbox2=\hbox{$\m@th#1\vcenter{\copy0}$}%
    \dimen0=\ht0
    \advance\dimen0 by-\ht2
    \ifdim\dimen0>0pt
        \mathinner{\raise\dimen0\hbox{$\m@th#1\left\lVert\vcenter{\copy0}\right\rVert$}}%
    \else
        \left\lVert#2\right\rVert
    \fi
    \endgroup
}
\newcommand{\norm}[1]{\mathinner{\mathpalette{\norm@auto}{#1}}}
\makeatother

\newcommand{\abs}[1]{\left | #1 \right |}
\newcommand{\D}{\mathbb{D}}
\newcommand{\C}{\mathbb{C}}
\newcommand{\R}{\mathbb{R}}
\newcommand{\dA}{\,\mathrm{d}A}
\newcommand{\dd}{\,\mathrm{d}}
\newcommand{\K}{\mathcal{K}}
\newcommand{\J}{\mathcal{J}}
\newcommand{\Jz}{\mathcal{J}_{0}}
\newcommand{\Lc}{\mathcal{L}}
\newcommand{\Hc}{\mathcal{H}}
\newcommand{\A}{\mathcal{A}}
\newcommand{\re}{\operatorname{Re}}
\newcommand{\im}{\operatorname{Im}}
\newcommand{\ol}{\overline}
\newcommand{\ee}{\mathrm{e}}

\title[The sharp H\"older constant in Mori's theorem]{The sharp H\"older constant\\ in Mori's theorem}
\author{Qiang Gao}
\address{School of Mathematical Sciences, Shenzhen University, 518060 Shenzhen, Guangdong, People's Republic of China;
    and Institute of Mathematics, Shenzhen University, 518060 Shenzhen, Guangdong, People's Republic of China}
\email{gaoqiangks@outlook.com}
\thanks{The author was supported by NSFC (No.~12401097) and the Guangdong Basic and Applied Basic Research Foundation (No.~2024A1515012362).}
\date{}

\begin{document}

\begin{abstract}
    Let $K\ge1$ and let $f$ be a $K$-quasiconformal homeomorphism of the unit disk $\D$ onto itself with $f(0)=0$.
    We prove that $\abs{f(a)-f(b)}\le16^{1-\frac{1}{K}}\abs{a-b}^{\frac{1}{K}}$ for all $a,b\in\D$.
    The constant is sharp, confirming the conjectured value of the Mori constant.
    The proof reduces the estimate along a Beltrami path to an infinitesimal inequality,
    which follows from quasiconformal distortion bounds and a three-pole integral estimate.
\end{abstract}

\maketitle

\setcounter{tocdepth}{1}
\tableofcontents

\section{Introduction}\label{sec:intro}

\subsection{Main result}\label{subsec:result}
Let $ \D = \set {z \in \C : \abs {z}<1}$. For $K \ge 1$ let $ \mathcal {F}_K$ be the family of orientation-preserving
$K$-quasiconformal homeomorphisms of $ \D $ onto itself with $f(0)=0$; equivalently, the
complex dilatation $ \mu _f=\frac{f_{ \bar {z}}}{f_z}$ satisfies $ \norm { \mu _f}_ \infty \le \frac{K-1}{K+1}$.
The \emph{Mori constant} is
\begin{equation}
    \label{eq:M}
    M(K)= \sup \set { \frac { \abs {f(a)-f(b)}}{ \abs {a-b}^{\frac{1}{K}}}: \  f \in \mathcal {F}_K, \  a,b \in \D , \  a \neq b}.
\end{equation}
Mori \cite{Mori} proved that $M(K) \le 16$, and $16$ is the best constant independent
of~$K$. For each fixed $K$, a classical extremal construction \cite[Ch.~II, \S3]{LV} gives
$M(K) \ge 16^{1-\frac{1}{K}}$. It has long been conjectured that this lower bound
is the exact value of $M(K)$ \cite[p.~68]{LV}.
In this paper, we give an affirmative answer to this conjecture.

\begin{theorem}[Main theorem]
    \label{thm:main}
    Let $K \ge 1$ and $f \in \mathcal {F}_K$. Then
    \begin{equation}
        \label{eq:main}
        \abs {f(a)-f(b)} \le 16^{1-\frac{1}{K}} \, \abs {a-b}^{\frac{1}{K}} \qquad (a,b \in \overline{\D} ).
    \end{equation}
    Consequently $M(K)=16^{1-\frac{1}{K}}$ for every $K \ge 1$.
\end{theorem}
\subsection{Previous results}
Semenov \cite[Thm.~9]{Sem87} used quasiconformal flows to obtain
$M(K)\le124^{1-\frac{1}{K}}$.
Fehlmann and Vuorinen \cite[Thm.~1.3]{FV} proved Hölder estimates whose coefficients tend to $1$ as $K \to 1$, so $M(K) \to 1$ (in every dimension).
Semenov \cite[Thm.~1]{Sem91} subsequently gave a deformation proof of the sharp radial estimate
$\abs{f(z)}\le4^{1-\frac{1}{K}}\abs{z}^{\frac{1}{K}}$, together with distortion estimates for $1-\abs{f(z)}^2$.
Qiu \cite[Thms.~1 and~2]{Qiu} proved
\begin{equation}
    M(K) \le 16^{1-\frac{1}{K}}R_Q(K), \qquad R_Q(K)= \min \set { \bigl ( {23}/{8} \bigr )^{1-\frac{1}{K}}, \  \bigl (1+2^{3-2K} \bigr )^{\frac{1}{K}}},
\end{equation}
which yields $M(K) \le 46^{1-\frac{1}{K}}$. Qiu also obtained the conjectured coefficient for pairs of points in a disk of radius $\sin 25.7^ \circ$.
Anderson, Vamanamurthy and Vuorinen \cite[Thm.~15.5]{AVV}
obtained coefficients for special configurations. Bhayo and Vuorinen \cite{BV} studied the analogous
problem in $ \R ^n$; see also the survey \cite[\S4]{Bh}. Semenov \cite[Question~2]{Sem} discussed a
variational formulation of the problem. Table~\ref{tab:comparison} summarises the planar results.
\begin{table}[ht]
    \centering\small
    \caption{Coefficients $C$ in $ \abs {f(a)-f(b)} \le C \abs {a-b}^{\frac{1}{K}}$ for $f \in \mathcal {F}_K$.}
    \label{tab:comparison}
    \begin{tabular}{lll}
        \toprule
        Result                       & Point pairs                                  & Coefficient $C$            \\
        \midrule
        Mori \cite{Mori}             & all $a,b \in \D $                            & $16$                       \\
        Semenov \cite[Thm.~9]{Sem87} & all $a,b\in\D$                               & $124^{1-\frac{1}{K}}$      \\
        Semenov \cite[Thm.~1]{Sem91} & $b=0$                                        & $4^{1-\frac{1}{K}}$        \\
        Qiu \cite[Thm.~1]{Qiu}       & all $a,b \in \D $                            & $16^{1-\frac{1}{K}}R_Q(K)$ \\
        Qiu \cite[Thm.~2]{Qiu}       & $ \abs {a}, \abs {b} \le \sin 25.7^ \circ $  & $16^{1-\frac{1}{K}}$       \\
        AVV \cite[Thm.~15.5(2)]{AVV} & $a,b \in \partial \D $                       & $16^{1-\frac{1}{K}}$       \\
        AVV \cite[Thm.~15.5(3)]{AVV} & $a \in \partial \D $ or $b \in \partial \D $ & $40^{1-\frac{1}{K}}$       \\
        AVV \cite[Thm.~15.5(4)]{AVV} & $0 \in [a,b]$                                & $8^{1-\frac{1}{K}}$        \\
        Theorem~\ref{thm:main}       & all $a,b \in \overline { \D }$               & $16^{1-\frac{1}{K}}$       \\
        \bottomrule
    \end{tabular}
\end{table}

\subsection{Strategy of the proof}\label{subsec:strategy}
Our proof strategy is based on the method of quasiconformal flow used by Semenov in \cite{Sem87, Sem91}.

\smallskip\noindent\emph{Reduction to an infinitesimal bound (Section~\ref{sec:reduction}).}
Fix $\mu\in L^{\infty}(\D)$ with $\norm{\mu}_{\infty}\le1$.
For $0\le t<1$, let $f_t:\D\to\D$ be the unique quasiconformal homeomorphism with complex dilatation $t\mu$
that fixes $0$ and $1$.
Put $K_t=\frac{1+t}{1-t}$. Then $f_0=\mathrm{id}$ and $f_t$ is $K_t$-quasiconformal.
Fix distinct $a,b\in\D$, and write
$\Lc_{\mu}(a,b):=\left.\frac{\dd}{\dd t}\log\abs{f_t(a)-f_t(b)}\right|_{t=0^+}$
for the initial rate of change of the logarithmic distance.
If \eqref{eq:main} holds for all admissible mappings, then
\begin{equation}
    \label{eq:main-converse-log}
    \log\abs{f_t(a)-f_t(b)}
    \le\left(1-\frac{1}{K_t}\right)\log16+\frac{1}{K_t}\log\abs{a-b}.
\end{equation}
Both sides are equal at $t=0$. Comparing their right derivatives there
yields
\begin{equation}
    \label{eq:main-converse-derivative}
    \Lc_{\mu}(a,b)\le2\log16-2\log\abs{a-b}
    =2\log\frac{16}{\abs{a-b}}=\Hc(a,b).
\end{equation}
This motivates the infinitesimal bound proved in Theorem~\ref{thm:kernel}.
Conversely, Section~\ref{subsec:integration} derives \eqref{eq:main} by applying this bound along a Beltrami path and integrating.
Thus the infinitesimal bound is equivalent to the main estimate for all $K\ge1$.

Semenov \cite[\S3, Lem.~4, pp.~203--204]{Sem87} used an integral representation of quasiconformal deformation fields
to derive a two-point estimate. Combining this estimate with quasiconformal flows, he obtained
$M(K)\le124^{1-\frac{1}{K}}$ \cite[Thm.~9, p.~205]{Sem87}.
In later work \cite[Question~2, p.~249]{Sem}, he explicitly formulated the associated integral problem,
observing that its sharp estimate would settle Mori's conjecture. In our notation, his integral is
\begin{equation}
    \label{eq:semenov-integral}
    I_{\mathrm{S}}(a,b)=\frac{r^2}{2\pi}\int_{\D}
    \abs{\frac{1}{(w-a)(w-b)}
        +\frac{\ol{a}+\ol{b}-\ol{ab}w}{w(1-\ol{a}w)(1-\ol{b}w)}}\dA(w).
\end{equation}
The connection with our infinitesimal bound follows from the integral representation \eqref{eq:L-disk-integral}:
\begin{equation}
    \label{eq:semenov-equivalence}
    I_{\mathrm{S}}(a,b)=\frac{r^2}{2}
    \sup_{\norm{\mu}_{\infty}\le1}\Lc_{\mu}(a,b).
\end{equation}

\smallskip\noindent\emph{The bound $\Lc\le\Hc$ (Sections~\ref{sec:linear}--\ref{sec:regions}).}
We estimate $\Lc=\Lc_{\mu}(a,b)$ using four bounds for the initial velocity $v$
(Proposition~\ref{prop:four}):
\begin{itemize}[leftmargin=2.2em]
    \item[(R)] the infinitesimal quasiconformal Schwarz lemma at $a$ and at $b$;
    \item[(S)] the infinitesimal two-point (hyperbolic) Schwarz lemma for the pair $(a,b)$;
    \item[(A)] the infinitesimal Agard--Gehring angular distortion bound for the triangle $(0,a,b)$;
    \item[(W)] a bound on the modulus of the complex quantity $W=\frac{v(a)}{a}-\frac{v(b)}{b}$ in terms of
          a three-pole integral $\J(a,b)$.
\end{itemize}
We bound the normalised integral $\sqrt{\abs{ab}}\,\J(a,b)$ by $\Jz(s)$, the corresponding integral whose two nonzero
poles lie on the unit circle, where $s:=\frac{\abs{a-b}}{\abs{a}+\abs{b}}$ is the relative chord length.
An explicit logarithmic series expansion of $\Jz(s)$ gives the estimates needed below (Section~\ref{sec:threepole}).
Combining (S) with (W) gives the first bound for $\Lc$, while combining (R) with (W) gives a second bound
under the geometric condition stated in Proposition~\ref{prop:joint2} (Section~\ref{sec:joint}).

\begin{figure}[ht]
    \centering
    \begin{tikzpicture}[
            box/.style={draw,rounded corners=2pt,align=center,font=\footnotesize,inner sep=4pt,text width=4.6cm},
            inp/.style={draw,dashed,rounded corners=2pt,align=center,font=\footnotesize,inner sep=4pt,text width=4.8cm,fill=black!3},
            >=Stealth, node distance=5mm and 8mm]
        \node[box] (vel) {Initial velocity\\ $v(z)=\left.\frac{\dd}{\dd t}f_t(z)\right|_{t=0^+}$\\ (Prop.~\ref{prop:velocity})};
        \node[box,below=of vel] (four) {Four estimates for $v$\\ (R), (S), (A), (W) \quad(Prop.~\ref{prop:four})};
        \node[box,below=of four] (joint) {First bound: (S) + (W)\\ Second bound: (R) + (W)\\ (Props.~\ref{prop:joint1}, \ref{prop:joint2})};
        \node[box,below=of joint] (bound) {Four regions and the origin case\\ $\Lc\le\Hc$ for every $\norm{\mu}_{\infty}\le1$\\ (Thm.~\ref{thm:kernel})};
        \node[box,below=of bound] (integ) {Integration along $f_t$\\ (\S\ref{subsec:integration})};
        \node[box,below=of integ,fill=blue!6] (main) {Theorem~\ref{thm:main}};
        \node[inp,left=of vel] (mrmt) {Measurable Riemann mapping\\ theorem, variation formula (\S\ref{subsec:beltrami})};
        \node[inp,left=of four] (cls) {QC Schwarz lemmas,\\ Agard--Gehring theorem,\\ elliptic estimates (\S\ref{subsec:elliptic})};
        \node[inp,left=of joint] (three) {Three-pole integral: comparison\\ and explicit series (\S\ref{sec:threepole})};
        \node[inp,left=of bound] (scal) {Elementary inequalities\\ and the radial Schwarz estimate\\ (\S\ref{sec:regions}, App.~\ref{app:elementary})};
        \node[inp,left=of main] (sharp) {Classical lower bound\\ (\S\ref{subsec:sharpness})};
        \draw[->] (mrmt) -- (vel);
        \draw[->] (cls) -- (four);
        \draw[->] (three) -- (joint);
        \draw[->] (scal) -- (bound);
        \draw[->] (vel) -- (four);
        \draw[->] (four) -- (joint);
        \draw[->] (joint) -- node[right,font=\scriptsize]{Regions II--IV} (bound);
        \draw[->] (four.east) -- ++(7mm,0)
        -- node[right,font=\scriptsize,align=left]{Region I\\ (A) + (R)} ([xshift=7mm]bound.east) -- (bound.east);
        \draw[->] (bound) -- (integ);
        \draw[->] (integ) -- (main);
        \draw[->] (sharp) -- (main);
    \end{tikzpicture}
    \caption{Structure of the proof.}
    \label{fig:flow}
\end{figure}
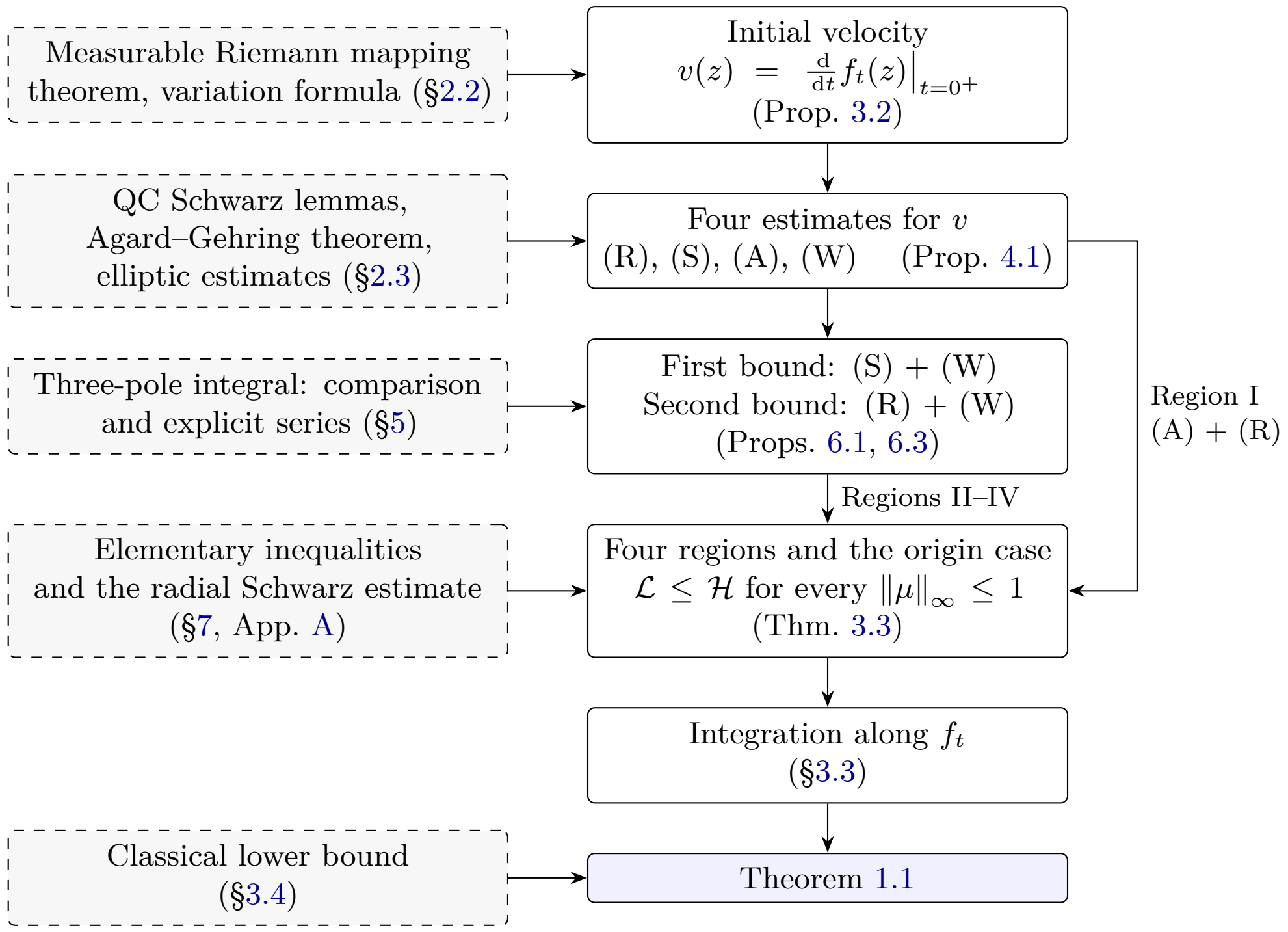

We divide the nonzero point pairs into four regions. Region~I uses (A) and (R) directly;
Regions~II and~III use the first bound; Region~IV uses the second bound and the monotonicity of its difference from $\Hc$ under
common scaling of $a$ and $b$. In each region, elementary inequalities in $\abs{a}$, $\abs{b}$, and $\abs{a-b}$
show that the chosen bound is at most $\Hc$ (Section~\ref{sec:regions} and Appendix~\ref{app:elementary}).
Pairs containing the origin are handled separately by the radial Schwarz estimate.
Figure~\ref{fig:flow} shows how the parts fit together.

\section{Preliminaries}\label{sec:classical}

In this section, we review several classical results and a number of lemmas. Most of the material presented here can be found in the textbooks \cite{Ahl,AVV,AIM,IT,LV,Mag}.

\subsection{Quasiconformal mappings}\label{subsec:qc}

Let $\Omega,\Omega'\subset\C$ be domains and let $z=x+\mathrm{i}y$.
For a map with weak first derivatives, we use the Wirtinger derivatives
\begin{equation}
    \label{eq:qc-wirtinger}
    f_z=\frac{1}{2}(f_x-\mathrm{i}f_y),\qquad
    f_{\bar{z}}=\frac{1}{2}(f_x+\mathrm{i}f_y).
\end{equation}
For $K\ge1$, an orientation-preserving homeomorphism $f:\Omega\to\Omega'$ is
$K$-quasiconformal if $f\in W^{1,2}_{\mathrm{loc}}(\Omega)$ and
\begin{equation}
    \label{eq:qc-definition}
    \abs{f_{\bar{z}}}\le\frac{K-1}{K+1}\abs{f_z}\quad\text{a.e. on }\Omega.
\end{equation}
Quasiconformal maps are differentiable a.e., and at almost every point, the differential $df$ maps circles to ellipses whose ratio of major to minor axes is at most $K$. A $1$-quasiconformal map is conformal. For a quasiconformal map, $f_z\ne0$ a.e., its \emph{Beltrami coefficient}, or
\emph{complex dilatation}, is defined a.e.\ by
\begin{equation}
    \label{eq:qc-beltrami}
    \mu_f=\frac{f_{\bar{z}}}{f_z},\qquad
    \norm{\mu_f}_{\infty}\le\frac{K-1}{K+1}<1.
\end{equation}

Let $\mu_t\in L^{\infty}(\Omega)$ be a family depending on a real or complex parameter $t$ near $t_0$.
We say that this family is \emph{differentiable at $t_0$ in the $L^{\infty}$ norm}, with derivative $\dot{\mu}\in L^{\infty}(\Omega)$, if
\begin{equation}
    \label{eq:linfty-derivative}
    \norm{\mu_{t_0+h}-\mu_{t_0}-h\dot{\mu}}_{\infty}=o(\abs{h})
    \quad\text{as }h\to0.
\end{equation}

The inverse of a $K$-quasiconformal map is $K$-quasiconformal, and the composition of a $K_1$-quasiconformal
map with a $K_2$-quasiconformal map is $K_1K_2$-quasiconformal.
Quasiconformality is invariant under conformal changes of coordinates; hence it is
also defined on the Riemann sphere $\widehat{\C}=\C\cup\set{\infty}$ using local conformal coordinates.
For quasiconformal maps $f:\Omega\to\Omega'$ and $g:\Omega\to\Omega''$, there holds
\begin{equation}
    \label{eq:qc-composition}
    \mu_{g\circ f^{-1}}(f(z))
    =\frac{\mu_g(z)-\mu_f(z)}{1-\ol{\mu_f(z)}\mu_g(z)}
    \frac{\,f_z(z)\,}{\ol{f_z(z)}}\quad\text{a.e. on }\Omega.
\end{equation}
In particular, if $\mu_f=\mu_g$ a.e., then $g\circ f^{-1}$ is conformal. %

The Jacobians of quasiconformal maps are positive a.e., and both the maps and their inverses preserve null sets. In complex notation the Jacobian is
\begin{equation}
    \label{eq:qc-jacobian}
    J_f=\abs{f_z}^2-\abs{f_{\bar{z}}}^2
    =(1-\abs{\mu_f}^2)\abs{f_z}^2>0\quad\text{a.e.}
\end{equation}
For a quasiconformal homeomorphism $f:\Omega\to\Omega'$, the area change of variables is
\begin{equation}
    \label{eq:qc-area}
    \int_{\Omega'}\Phi(w)\dA(w)
    =\int_{\Omega}\Phi(f(z))J_f(z)\dA(z),\qquad \Phi\in L^1(\Omega').
\end{equation}

Finally, every quasiconformal homeomorphism of $\D$ onto itself extends to a homeomorphism of
$\overline{\D}$ onto itself. This allows the estimates later in this paper to be passed between interior and boundary points by continuity.

\subsection{The Beltrami equation}\label{subsec:beltrami}

\begin{theorem}[{\cite[Thm.~4.30, p.~102]{IT}}]
    \label{thm:mrmt}
    For every $\mu\in L^{\infty}(\C)$ with $\norm{\mu}_{\infty}<1$, there is a homeomorphism $f:\widehat{\C}\to\widehat{\C}$ which is a quasiconformal mapping of $\mathbb{C}$ with complex dilatation $\mu$. Moreover, $f$ is uniquely determined by the following normalization conditions:
    \begin{equation}
        \label{eq:normalisation}
        f(0)=0,\qquad f(1)=1,\qquad f(\infty)=\infty.
    \end{equation}
\end{theorem}
Following \cite{IT}, we denote this normalised map by $f^{\mu}$.
\begin{theorem}[{\cite[Cor.\ to Thm.~4.37, p.~110]{IT}}]
    \label{thm:variation}
    Let $\mu_h$ be a family of Beltrami coefficients depending on a real or complex parameter $h$ near $0$,
    with $\mu_h\in L^{\infty}(\C)$ and $\norm{\mu_h}_{\infty}<1$. Suppose that
    \begin{equation}
        \label{eq:variation-coefficient}
        \mu_h=\mu+h\nu+h\varepsilon_h,
        \qquad \norm{\varepsilon_h}_{\infty}\longrightarrow0
        \quad\text{as }h\to0,
    \end{equation}
    where $\mu,\nu,\varepsilon_h\in L^{\infty}(\C)$ and $\norm{\mu}_{\infty}<1$.
    Then
    \begin{equation}
        \label{eq:variation-original}
        f^{\mu_h}(z)=f^{\mu}(z)+h\dot{f}^{\mu}[\nu](z)+o(\abs{h})
        \quad\text{as }h\to0,
    \end{equation}
    locally uniformly for $z\in\C$, where, for $z\in\C\setminus\set{0,1}$,
    \begin{equation}
        \label{eq:variation-original-integral}
        \dot{f}^{\mu}[\nu](z)
        =-\frac{1}{\pi}\int_{\C}\nu(\zeta)
        \frac{f^{\mu}(z)(f^{\mu}(z)-1)((f^{\mu})_{\zeta}(\zeta))^2}
        {f^{\mu}(\zeta)(f^{\mu}(\zeta)-1)(f^{\mu}(\zeta)-f^{\mu}(z))}\dA(\zeta),
    \end{equation}
    and $\dot{f}^{\mu}[\nu](0)=\dot{f}^{\mu}[\nu](1)=0$.
\end{theorem}

We will use the following form of \eqref{eq:variation-original-integral}, with the integral taken in the image plane.
\begin{corollary}[Variation formula in the image plane]
    \label{cor:transported-variation}
    Let $\mu_h$, $\mu$, and $\nu$ satisfy the assumptions of Theorem~\ref{thm:variation}.
    Define $\dot{\lambda}$ a.e.\ by
    \begin{equation}
        \label{eq:transported}
        \dot{\lambda}(f^{\mu}(z))
        =\frac{\nu(z)}{1-\abs{\mu(z)}^2}\cdot\frac{\, (f^{\mu})_z(z) \,}{\ol{(f^{\mu})_z(z)}}.
    \end{equation}
    Then
    \begin{equation}
        \label{eq:transported-expansion}
        f^{\mu_h}(z)=f^{\mu}(z)+hV[\dot{\lambda}](f^{\mu}(z))+o(\abs{h})
        \quad\text{as }h\to0,
    \end{equation}
    locally uniformly for $z\in\C$. Here $V$ is the operator on $L^{\infty}(\C)$ defined, for $y\in\C\setminus\set{0,1}$, by
    \begin{equation}
        \label{eq:V}
        V[\nu](y)=-\frac{1}{\pi}\int_{\C}\nu(w)\frac{y(y-1)}{w(w-1)(w-y)}\dA(w),
    \end{equation}
    with $V[\nu](0)=V[\nu](1)=0$.
\end{corollary}

\begin{proof}
    Since $f^{\mu}$ and its inverse preserve null sets, \eqref{eq:transported} defines a function $\dot{\lambda}\in L^{\infty}(\C)$, with
    \begin{equation}
        \label{eq:image-variation-norm}
        \norm{\dot{\lambda}}_{\infty}
        \le\frac{\norm{\nu}_{\infty}}{1-\norm{\mu}_{\infty}^2}.
    \end{equation}
    By the Jacobian identity \eqref{eq:qc-jacobian},
    \begin{equation}
        \label{eq:variation-coefficient-change}
        \begin{aligned}
            & \frac{\nu(\zeta)((f^{\mu})_{\zeta}(\zeta))^2}{J_{f^{\mu}}(\zeta)} \\
            ={} & \frac{\nu(\zeta)}{1-\abs{\mu(\zeta)}^2}
            \frac{\,(f^{\mu})_{\zeta}(\zeta)\,}{\ol{(f^{\mu})_{\zeta}(\zeta)}}
            =\dot{\lambda}(f^{\mu}(\zeta)).
        \end{aligned}
    \end{equation}
    The change of variables $w=f^{\mu}(\zeta)$ in \eqref{eq:variation-original-integral}, using \eqref{eq:qc-area},
    therefore gives $\dot{f}^{\mu}[\nu](z)=V[\dot{\lambda}](f^{\mu}(z))$.
    The expansion \eqref{eq:transported-expansion}, including its local uniformity, now follows from \eqref{eq:variation-original}.
\end{proof}

\subsection{Some functions}\label{subsec:elliptic}\label{app:geometry}\label{subsec:distortion}
We use the complete elliptic integral of the first kind $\K$, the Gr\"{o}tzsch modulus $\mathcal{G}$,
and the distortion function $\varphi_K$; see \cite[(3.1), p.~48; (5.1), p.~80; (10.1), p.~203]{AVV}.
For $0<q<1$ put $q'= \sqrt {1-q^2}$ and
\begin{equation}
    \label{eq:elliptic}
    \begin{aligned}
        \K (q)          & = \int _0^{ \frac{\pi}{2}} \frac { \dd \theta }{ \sqrt {1-q^2 \sin ^2 \theta }}, \\
        \mathcal{G} (q) & = \frac { \pi }{2} \, \frac { \K (q')}{ \K (q)},                                 \\
        \varphi _K(q)   & = \mathcal{G} ^{-1} \Bigl ( \frac { \mathcal{G} (q)}{K} \Bigr ) \quad (K>0).
    \end{aligned}
\end{equation}
The Grötzsch modulus function is generally denoted by $\mu$ in the literature. However, in this paper we denote it by $\mathcal{G}$ in order to reserve $\mu$ for the Beltrami differential.
The function $ \mathcal{G} $ is a strictly decreasing homeomorphism from $(0,1)$ onto $(0, \infty )$, and $\varphi_K$ is a strictly increasing homeomorphism from $(0,1)$ onto itself.
We also use
\begin{equation}
    \label{eq:Bm}
    \begin{aligned}
        B(q) & = \frac {4}{ \pi } \K (q) \K (q'),                                      \\
        m(q) & = \frac {q'^2}{2}B(q)= \frac {2}{ \pi } \, q'^2 \K (q) \K (q'),         \\
        c_0  & = \frac {1}{8} \exp \Bigl ( \frac {B(\frac{1}{ \sqrt {2}})}{2} \Bigr ).
    \end{aligned}
\end{equation}
Note that $B(q)=B(q')$, and $c_0= \frac {1}{8} \exp ( \frac {2} {\pi} \K (\frac{1}{ \sqrt {2}})^2) \approx 1.11516$.

We now give several geometric interpretations of these functions.

\emph{The elliptic function and Rectangles.} The Schwarz--Christoffel map
$z \mapsto \frac {1}{2} \int _0^z \frac { \dd u}{ \sqrt {u(1-u)(1-q^2u)}}$ sends the upper half-plane onto the rectangle
$(0, \K (q)) \times (0, \K (q'))$, with the boundary points $0,1,\frac{1}{q^2}, \infty $
mapped to its four vertices \cite[Ch.~6, \S2.3, pp.~238--239]{AhlCA}.
The square root is chosen analytic in the upper half-plane with positive boundary values on $(0,1)$.
The rectangle's height-to-width ratio is
$\frac{2 \mathcal{G} (q)}{\pi} $, and its area is $ \frac { \pi }{4}B(q)$.

\emph{The functions $\mathcal{G}$ and $ \varphi _K$.}
The Gr\"{o}tzsch modulus $\mathcal{G}(q)$ is the modulus of the slit disk $\D\setminus[0,q]$. By Gr\"{o}tzsch's extremal theorem, this is the greatest modulus among all ring domains contained in $\D$
and separating $0$ and $q$ from $\partial\D$ \cite[Ch.~II, \S1.1, pp.~53--54]{LV}.
\par The function $\varphi_K$ is a radial distortion bound for $K$-quasiconformal self-maps of the unit disk fixing the origin \cite[Ch.~II, \S3.1, pp.~63--64]{LV}. Let $f:\D\to\D$ be $K$-quasiconformal with $K\ge1$ and $f(0)=0$, and fix $z\in\D\setminus\set{0}$.
Slitting $\D$ along the segment $[0,z]$ gives a ring of modulus $\mathcal{G}(\abs{z})$.
By Gr\"{o}tzsch's extremal theorem, its image has modulus
at most $\mathcal{G}(\abs{f(z)})$. Combining this with the distortion of ring moduli under a $K$-quasiconformal map gives%
\begin{equation}
    \label{eq:radial-modulus-origin}
    \mathcal{G}(\abs{z})\le K\mathcal{G}(\abs{f(z)}).
\end{equation}
Since $\mathcal{G}$ is strictly decreasing, this is equivalent to
\begin{equation}
    \label{eq:radial-distortion-origin}
    \abs{f(z)}\le\mathcal{G}^{-1}\left(\frac{\mathcal{G}(\abs{z})}{K}\right)
    =\varphi_K(\abs{z}),
\end{equation}
which is the radial conclusion of Theorem~\ref{thm:schwarz}.

\emph{The functions $m$ and $B$.} As recorded in Lemma~\ref{lem:elliptic}(a) below (see also \cite[p.~37, (16)]{Qiu}), $ \varphi _{1+ \varepsilon }(q)=q(1+ \varepsilon \, m(q)+o( \varepsilon ))$ for fixed $q$ as $\varepsilon\downarrow0$.
Thus $m(q)$ describes the first-order relative change of the radial distortion bound $\varphi_K(q)$
as $K$ increases from $1$. The function $B$ enters the infinitesimal estimates obtained by differentiating the quasiconformal Schwarz lemma.
Both applications are developed in Section~\ref{sec:linear}.

\begin{figure}[ht]
    \centering
    \includegraphics[width=0.54\textwidth]{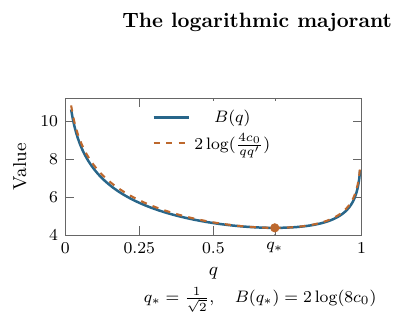}
    \caption{$B(q)$ and its logarithmic majorant $2 \log (\frac{4c_0}{qq'})$,
        which meet at $q_*=\frac{1}{ \sqrt {2}}$.}
    \label{fig:majorant}
\end{figure}

\emph{The constant $c_0$.} The majorant $2 \log (\frac{4c_0}{qq'})$ of $B(q)$ equals $B$ at the point $q=\frac{1}{ \sqrt {2}}$ (Figure~\ref{fig:majorant}).


We state two distortion theorems which play a central role in this paper.
\begin{theorem}[Quasiconformal Schwarz lemma {\cite[(3.9)~p.65]{LV}} {\cite[(10.18), p.~210]{AVV}}]
    \label{thm:schwarz}
    Let $g: \D \to \D $ be $K$-quasiconformal \textup{(}not necessarily onto\textup{)}. Then for $a,b \in \D $,
    \begin{equation}
        \delta (g(a),g(b)) \le \varphi _K( \delta (a,b)), \quad \text{where} \quad \delta (a,b)= \frac{\abs {a-b}}{ \abs {1- \ol {a}b}}.
    \end{equation}
    In particular, if $g(0)=0$ then
    \begin{equation}
        \abs {g(z)} \le \varphi _K( \abs {z}).
    \end{equation}
\end{theorem}

\begin{theorem}[Angular distortion; Agard--Gehring {\cite[Thm.~2]{AG}}{\cite[Thm.~14.22]{AVV}}]
    \label{thm:angle}
    Let $G$ be a $K$-quasiconformal homeomorphism of the Riemann sphere with $G( \infty )= \infty $. For distinct
    finite points $a_0,a_1,a_2$ define $ \Theta \in (0, \pi ]$ by
    \begin{equation}
        \sin \frac {\Theta} {2}= \frac { \abs {a_1-a_2}}{ \abs {a_0-a_1}+ \abs {a_0-a_2}},
    \end{equation}
    and define $ \widetilde {\Theta} $ in the same way for $G(a_0),G(a_1),G(a_2)$. Then
    $ \sin ( \frac{\widetilde {\Theta}}{4}) \le \varphi _K \bigl ( \sin ( \frac{\Theta}{4}) \bigr )$.
\end{theorem}

Note that the function denoted by $\varphi_K$ in \cite[(4.16), p.~10]{AG} is the inverse of our $\varphi_K$.


The following lemma forms the basis of the infinitesimal estimates developed later in the paper.
\begin{lemma}
    \label{lem:elliptic}
    \label{lem:c0}
    For $0<q<1$:
    \begin{enumerate}[label=\textup{(\alph*)}]
        \item $ \displaystyle \frac { \partial }{ \partial K} \log \varphi _K(q) \Big |_{K=1}=m(q)$;
        \item $m(q)+ \log q$ is a strictly decreasing function from $(0, 1)$ onto $(0, \log4)$,
              in particular $m(q) \le \log (\frac{4}{q})$;
        \item $B(q) \le 2 \log \frac {4c_0}{qq'}$, with equality at $q=\frac{1}{ \sqrt {2}}$;
        \item $ \log c_0<\frac{7}{60}$ and $c_0<\frac{9}{8}$.
    \end{enumerate}
\end{lemma}

\begin{proof}
    (a) By \cite[(10.6), p.~205]{AVV}, for fixed $0<q<1$ and $K>0$,
    \begin{equation}
        \label{eq:avv-phi-KK}
        \frac{\partial\varphi_K(q)}{\partial K}
        =\frac{4}{\pi^2K^2}ss'^2\K(s)^2\mu(q),
    \end{equation}
    where $s=\varphi_K(q)$, $s'=\sqrt{1-s^2}$ denotes the complementary modulus, and
    its function $\mu$ is our $\mathcal{G}$. Therefore, the identity \eqref{eq:avv-phi-KK} gives
    \begin{equation}
        \label{eq:avv-phi-K}
        \frac{\partial\varphi_K(q)}{\partial K}
        =\frac{4}{\pi^2K^2}\varphi_K(q)(1-\varphi_K^2(q))\K(\varphi_K(q))^2\mathcal{G}(q),
    \end{equation}
    Dividing \eqref{eq:avv-phi-K} by $\varphi_K(q)$ and setting $K=1$, so that $s=\varphi_1(q)=q$, gives
    \begin{equation}
        \left.\frac{\partial}{\partial K}\log\varphi_K(q)\right|_{K=1} = \frac{\left.\partial_K \varphi_K(q)\right|_{K=1}}{q}
        =\frac{4}{\pi^2}q'^2\K(q)^2\mathcal{G}(q)
        =\frac{2}{\pi}q'^2\K(q)\K(q')=m(q),
    \end{equation}
    where the last two equalities follow from \eqref{eq:elliptic} and \eqref{eq:Bm}.

    (b) This follows from \cite[Thm.~3.30(1), (6), p.~61]{AVV}.

    (c) The distortion estimate of Qiu and Vuorinen \cite[Thm.~1.3]{QV96} states that, for $t>0$ and $K \ge 1$,
    \begin{equation}
        \Bigl [ \frac { \varphi _K(q)}{ \varphi _{\frac{1}{K}}(q')} \Bigr ]^2 \le 16^{K-1}c_0^{2(1-\frac{1}{K})}t^{\frac{1}{K}}(1+t)^{K-\frac{1}{K}},
        \qquad q= \sqrt {\frac{t}{1+t}}.
    \end{equation}
    Both sides are equal at $K=1$, so the derivative of the left-hand side in $K$ at $K=1$ is no greater than that of the right-hand side. By (a), the
    derivative of the logarithm of the left side is $2m(q)+2m(q')=q'^2B(q)+q^2B(q')=B(q)$. Since $t=\frac{q^2}{q'^2}$
    and $1+t=\frac{1}{q'^2}$, the derivative of the logarithm of the right side is
    $ \log 16+2 \log c_0- \log t+2 \log (1+t)=2 \log (\frac{4c_0}{qq'})$. Equality at $q=\frac{1}{ \sqrt {2}}$ is the definition of $c_0$.

    (d) By the series expansion \eqref{eq:K-series},
    $\frac{2}{\pi}\K(\frac{1}{\sqrt{2}})=\sum_{n\ge0}a_n$, where
    $a_n=\frac{1}{2^n}\bigl(\frac{\binom{2n}{n}}{4^n}\bigr)^2$.
    Since $\frac{a_{n+1}}{a_n}<\frac{1}{2}$, the tail starting at $a_4$ is less than $2a_4$, so
    \begin{equation}
        \frac{2}{\pi}\K(\frac{1}{\sqrt{2}})<1+\frac{1}{8}+\frac{9}{256}+\frac{25}{2048}+\frac{1225}{131072}<\frac{591}{500}.
    \end{equation}
    Using $\log2>\frac{842}{1215}$ and $\pi<\frac{22}{7}$, we obtain
    \begin{equation}
        \log c_0=\frac{\pi}{2}\Bigl(\frac{2}{\pi}\K(\frac{1}{\sqrt{2}})\Bigr)^2-3\log2
        <\frac{11}{7}\Bigl(\frac{591}{500}\Bigr)^2-\frac{842}{405}<\frac{7}{60}.
    \end{equation}
    Finally, $\frac{7}{60}<\log\frac{9}{8}$, and hence $c_0<\frac{9}{8}$.
\end{proof}

\subsection{Series and integrals}\label{subsec:series}
The rising factorial (Pochhammer symbol) is defined by
\begin{equation}
    ( \alpha )_0=1 \quad\text{and}\quad
    ( \alpha )_n= \alpha ( \alpha +1) \cdots ( \alpha +n-1)
    \quad \text{for integers}\quad n \ge 1.
\end{equation}
Put
\begin{equation}
    \label{eq:bn}
    b_n= \binom {2n}{n}4^{-n}= \frac {(\frac{1}{2})_n}{n!}, \qquad d_n=H_{2n}-H_n, \qquad H_n= \sum _{j=1}^n \frac {1}{j}
    \quad (d_0=0).
\end{equation}
For $n \ge 1$, $0<d_n= \sum _{j=n+1}^{2n} \frac {1}{j}< \int _n^{2n} \frac { \dd u}{u}= \log 2$, and $d_n$ increases
with $n$ (because $d_{n+1}-d_n= \frac {1}{2n+1}- \frac {1}{2n+2}$), with $d_1= \frac {1}{2}$. We use the
following classical formulas \cite[8.113.1, 8.113.3, 3.621.3, 4.387.4, 1.110, 1.513.1]{GR}:
\begin{equation}
    \label{eq:K-series}
    \K (q)= \frac {\pi} {2} \sum _{n \ge 0}b_n^2q^{2n} \qquad (0 \le q<1),
\end{equation}
\begin{equation}
    \label{eq:compl-series}
    \K \bigl ( \sqrt {1- q ^2} \bigr )= \sum _{n \ge 0}b_n^2 q ^{2n} \Bigl ( \log \frac {4} {q} -2d_n \Bigr )
    \qquad (0< q <1),
\end{equation}
\begin{equation}
    \label{eq:sine-moments}
    \int _0^{ \frac{\pi}{2}} \! \sin ^{2n} \theta \dd \theta = \frac {\pi} {2} b_n, \qquad
    \int _0^{ \frac{\pi}{2}} \! \sin ^{2n} \theta \log \sin \theta \dd \theta = \frac {\pi} {2} b_n(d_n- \log 2),
\end{equation}
\begin{equation}
    \label{eq:binomial}
    (1-y)^{- \alpha }= \sum _{n \ge 0} \frac {( \alpha )_n}{n!}y^n \quad ( \alpha >0, \  \abs {y}<1),
\end{equation}
\begin{equation}
    \label{eq:log-series}
    \log \frac {1+y}{1-y}=2 \sum _{n \ge 0} \frac {y^{2n+1}}{2n+1} \quad (0<y<1).
\end{equation}
For reference, the coefficient values used in this paper are
\begin{equation}
    \label{eq:b-explicit-values}
    b_0=1,\qquad b_1=\frac{1}{2},\qquad b_2=\frac{3}{8},\qquad
    b_3=\frac{5}{16},\qquad b_4=\frac{35}{128},
\end{equation}
\begin{equation}
    \label{eq:d-explicit-values}
    d_0=0,\qquad d_1=\frac{1}{2}.
\end{equation}

\section{The infinitesimal bound and the proof of the main theorem}\label{sec:reduction}

\subsection{Reflection symmetry}\label{subsec:reflection}
For $\mu\in L^{\infty}(\D)$, define
\begin{equation}
    \label{eq:reflection}
    \mu ^{ \mathrm {e}}(z)=
    \begin{cases}
        \mu (z), & \abs {z}<1, \\[2pt] \frac {z^2}{ \bar {z}^2} \, \ol { \mu (\frac{1}{\bar {z}})},& \abs {z}>1.
    \end{cases}
\end{equation}
For real $t$ with $ \abs {t} \norm { \mu }_ \infty <1$, let $F_t=f^{t\mu^{\mathrm{e}}}$ be the normalised solution
(Theorem~\ref{thm:mrmt}) with coefficient $t \mu ^{ \mathrm {e}}$. Let $\iota(z)=\frac{1}{\bar{z}}$ be reflection in the unit circle.

\begin{lemma}[Reflection symmetry]
    \label{lem:disk-preserving}
    The reflection extension and the normalised maps above have the following properties.
    \begin{enumerate}[label=(\roman*),ref=(\roman*)]
        \item\label{item:reflection-maps}
              Let $G$ be a quasiconformal homeomorphism of the sphere fixing $0$, $1$, and $\infty$. Then $G$ is reflection-symmetric if and only if $\mu_G$ is reflection-symmetric:
              \begin{equation}
                  \label{eq:reflection-equivalence}
                  G\circ\iota=\iota\circ G
                  \quad\Longleftrightarrow\quad
                  \mu_G=(\mu_G|_{\D})^{\mathrm{e}}\quad\text{a.e.}
              \end{equation}
              In this case, $G$ maps $\D$ onto itself.
              In particular, $F_t\circ\iota=\iota\circ F_t$, and $f_t:=F_t|_{\D}$ is a
              $\frac{1+\abs{t}\norm{\mu}_{\infty}}{1-\abs{t}\norm{\mu}_{\infty}}$-quasiconformal self-map of $\D$ fixing $0$.
        \item\label{item:reflection-variation}
              Suppose that the base coefficient $\mu$ and the variation direction $\nu$ in Corollary~\ref{cor:transported-variation} are both reflection-symmetric.
              Then the image-plane coefficient $\dot{\lambda}$ defined by \eqref{eq:transported} is also reflection-symmetric:
              \begin{equation}
                  \label{eq:reflection-variation}
                  \dot{\lambda}=(\dot{\lambda}|_{\D})^{\mathrm{e}}\quad\text{a.e. on }\C.
              \end{equation}
    \end{enumerate}
\end{lemma}

\begin{proof}
    \ref{item:reflection-maps} Put
    \begin{equation} \label{eq:def:H}
        H=\iota\circ G\circ\iota.
    \end{equation}
    As $\iota$ is anticonformal, $H$ is again an orientation-preserving quasiconformal homeomorphism. Direct computation shows that
    \begin{equation}
        (G\circ \iota)_z =  G_{\bar{z}} \bar{\iota}_z,\quad (G\circ \iota)_{\bar{z}} =  G_{{z}} \iota_{\bar{z}},
    \end{equation}
    and
    \begin{equation}
        \label{eq:reflection-derivatives}
        H_z = \iota_{\bar{z}} \cdot (\overline{G\circ \iota})_{z},\quad H_{\bar{z}} = \iota_{\bar{z}} \cdot (\overline{G\circ \iota})_{\bar{z}}.
    \end{equation}
    Therefore,
    \begin{equation}
        \label{eq:reflection-chain}
        \mu_H(z) = \frac{(\overline{G\circ \iota})_{\bar{z}}} {(\overline{G\circ \iota})_{{z}}} = \frac{\overline{G_{\bar{z}} \bar{\iota}_z}}{ \overline{G_z \iota_{\bar{z}}}} = \overline{\mu_G(\iota(z))} \frac{z^2}{\bar{z}^2}.
    \end{equation}
    \par Suppose that $G$ is reflection-symmetric. Then $H=G$, so $\mu_H=\mu_G$, which is exactly $\mu_G = (\mu_G|_{\D})^{\mathrm{e}}$.
    \par Conversely, suppose that $\mu_G = (\mu_G|_{\D})^{\mathrm{e}}$. Let $H$ be defined as in \eqref{eq:def:H}. Then $\mu_H = \mu_G$ by \eqref{eq:reflection-chain}. Since $H$ also fixes $0$, $1$, and $\infty$, $G \circ H^{-1}$ is the identity map. The conclusion follows.

    \ref{item:reflection-variation} By \ref{item:reflection-maps}, the normalised map $f^{\mu}$ satisfies
    $f^{\mu}(\iota(z))=\frac{1}{\ol{f^{\mu}(z)}}$.
    Differentiating with respect to $\bar{z}$ gives, a.e.\ for $z\ne0$,
    \begin{equation}
        \label{eq:path-reflected-derivative}
        (f^{\mu})_z(\iota(z))
        =\frac{\bar{z}^2}{\ol{f^{\mu}(z)}^2}\ol{(f^{\mu})_z(z)}.
    \end{equation}
    Since $\mu$ and $\nu$ are reflection-symmetric, we have
    $\abs{\mu(\iota(z))}=\abs{\mu(z)}$ and
    $\nu(\iota(z))=\frac{z^2}{\bar{z}^2}\ol{\nu(z)}$.
    Substituting these identities into \eqref{eq:transported}, and using $f^{\mu}\circ\iota=\iota\circ f^{\mu}$, we obtain
    \begin{equation}
        \label{eq:path-reflected-coefficient}
        \begin{aligned}
            & \dot{\lambda}(\iota(f^{\mu}(z))) \\
            ={} & \frac{f^{\mu}(z)^2}{\ol{f^{\mu}(z)}^2}
            \frac{\ol{\nu(z)}}{1-\abs{\mu(z)}^2}
            \frac{\ol{(f^{\mu})_z(z)}}{\,(f^{\mu})_z(z)\,} \\
            ={} & \frac{f^{\mu}(z)^2}{\ol{f^{\mu}(z)}^2}\ol{\dot{\lambda}(f^{\mu}(z))}.
        \end{aligned}
    \end{equation}
    Since $f^{\mu}$ preserves null sets, setting $w=f^{\mu}(z)$ gives the reflection identity for $\dot{\lambda}$ a.e.\ in the image plane.
    This proves \eqref{eq:reflection-variation}.
\end{proof}

\begin{proposition}
    \label{prop:velocity}
    Let $\mu\in L^{\infty}(\D)$ and $f_t$ be defined as in Lemma~\ref{lem:disk-preserving}. Locally uniformly on $\D$, $f_t(z)=z+t\,v(z)+o(t)$, where
    \begin{equation}
        \label{eq:velocity}
        v=V[\mu^{\mathrm{e}}]|_{\D}.
    \end{equation}
    Moreover, $v(0)=0$, and the velocity depends real-linearly on $\mu$.
\end{proposition}

\begin{proof}
    Apply Theorem~\ref{thm:variation} to the coefficient family $h\mapsto h\mu^{\mathrm{e}}$ at the zero coefficient.
    Since $f^{0}=\mathrm{id}$ and $V[\mu^{\mathrm{e}}]=\dot{f}^{0}[\mu^{\mathrm{e}}]$, this gives the stated locally uniform expansion and \eqref{eq:velocity}.
    The reflection extension is real-linear and $V$ is linear, so the velocity is real-linear in $\mu$.
    Finally, $v(0)=0$ follows from $f_t(0)=0$.
\end{proof}

\subsection{Infinitesimal bound}
For distinct $a,b \in \D $, define the chord log-speed and the target by
\begin{equation}
    \label{eq:L-kernel}
    \Lc _ \mu (a,b)= \re \frac {v(a)-v(b)}{a-b}, \qquad
    \Hc (a,b)=2 \log \frac {16}{ \abs {a-b}}.
\end{equation}
When the direction and the pair are fixed, we write $ \Lc $ and $ \Hc $. The functional $ \Lc _ \mu (a,b)$ is
real-linear in $ \mu $, and $ \Hc (a,b)>0$ because $0< \abs {a-b}<2$.
Splitting the integral in \eqref{eq:V} into its interior and exterior parts, using \eqref{eq:reflection}
and the substitution $w=\frac{1}{\ol{\zeta}}$ in the exterior part, and taking the real part of the divided difference gives
\begin{equation}
    \label{eq:L-disk-integral}
    \Lc_{\mu}(a,b)=-\frac{1}{\pi}\re\int_{\D}\mu(w)
    \left(\frac{1}{(w-a)(w-b)}
    +\frac{\ol{a}+\ol{b}-\ol{ab}w}{w(1-\ol{a}w)(1-\ol{b}w)}\right)\dA(w).
\end{equation}

\begin{theorem}[Infinitesimal bound]
    \label{thm:kernel}
    For every $ \mu \in L^ \infty ( \D )$ with $ \norm { \mu }_ \infty \le 1$ and all distinct $a,b \in \D $,
    \begin{equation}
        \label{eq:kernel}
        \Lc _ \mu (a,b) \le \Hc (a,b).
    \end{equation}
\end{theorem}

In the proofs of the two-point estimates, we take $a \ne b$ and label the points so that
$ \abs {a} \ge \abs {b}$, interchanging them if necessary.
In particular, $a \ne 0$.
The case $b=0$ in Theorem~\ref{thm:kernel} will be treated separately at the end of
Section~\ref{sec:regions}; the auxiliary estimates in Sections~\ref{sec:linear}--\ref{sec:regions}
are developed under the additional assumption $b \ne 0$.

Theorem~\ref{thm:kernel} is proved in Sections~\ref{sec:linear}--\ref{sec:regions}. We first show that it
implies Theorem~\ref{thm:main}.

\subsection{Proof of Theorem~\ref{thm:main} from Theorem~\ref{thm:kernel}}\label{subsec:integration}
If $a=b$ there is nothing to prove. If $K=1$, then $f$ is a conformal automorphism of $ \D $ fixing $0$, hence a
rotation by Schwarz's lemma, and \eqref{eq:main} is an equality. So let $a \ne b$ and $K>1$,
and label the points so that $ \abs {a} \ge \abs {b}$. Put
\begin{equation}
    \label{eq:normalized-beltrami-direction}
    k= \frac {K-1}{K+1} \in (0,1), \qquad \mu = \frac { \mu _f}{k}, \qquad \norm { \mu }_ \infty \le 1.
\end{equation}
Let $F_t$ ($0 \le t\le k$) be the normalised solutions with coefficient $t \mu ^{ \mathrm{e}}$, and $f_t=F_t|_ \D $.

For each fixed $t\in(0,k)$, apply Corollary~\ref{cor:transported-variation}
to the coefficient family $h\mapsto(t+h)\mu^{\mathrm{e}}$,
with base coefficient $t\mu^{\mathrm{e}}$ and variation direction $\mu^{\mathrm{e}}$.
The image-plane coefficient is given by
\begin{equation}
    \label{eq:path-image-coefficient}
    \dot{\lambda}(F_t(z))
    =\frac{\mu^{\mathrm{e}}(z)}{1-t^2\abs{\mu^{\mathrm{e}}(z)}^2}
    \frac{\,(F_t)_z(z)\,}{\ol{(F_t)_z(z)}},
\end{equation}
and satisfies
\begin{equation}
    \label{eq:path-image-norm}
    \norm{\dot{\lambda}}_{\infty}\le\frac{1}{1-t^2}.
\end{equation}
The resulting velocity formula is
\begin{equation}
    \label{eq:path-velocity}
    \frac{\dd}{\dd t}F_t(z)=V[\dot{\lambda}](F_t(z)).
\end{equation}

Since $t$ is real, both the base coefficient $t\mu^{\mathrm{e}}$ and the variation direction $\mu^{\mathrm{e}}$ are reflection-symmetric.
Lemma~\ref{lem:disk-preserving}\ref{item:reflection-variation} therefore gives
$\dot{\lambda}=\xi_t^{\mathrm{e}}$, where $\xi_t=\dot{\lambda}|_{\D}$.
By Proposition~\ref{prop:velocity}, $V[\dot{\lambda}]|_{\D}$ is the initial velocity field associated with the direction $\xi_t$.

The maps $f_t$ are injective and preserve $ \D $, so
$0< \abs {f_t(a)-f_t(b)}<2$ and $ \Hc (f_t(a),f_t(b))>0$.
The labels $a,b$ remain fixed along the path. By symmetry, Theorem~\ref{thm:kernel} applies to
$f_t(a),f_t(b)$ regardless of their relative moduli.
Set $\widehat{\xi}_t=(1-t^2)\xi_t$. The norm bound above gives $\norm{\widehat{\xi}_t}_{\infty}\le1$.
The velocity formula, Lemma~\ref{lem:log-modulus-derivative}, Theorem~\ref{thm:kernel}, and real-linearity yield
\begin{equation}
    \begin{aligned}
        \frac{\dd}{\dd t}\log\abs{f_t(a)-f_t(b)}
         & =\Lc_{\xi_t}\bigl(f_t(a),f_t(b)\bigr)                          \\
         & =\frac{1}{1-t^2}\Lc_{\widehat{\xi}_t}\bigl(f_t(a),f_t(b)\bigr) \\
         & \le\frac{2}{1-t^2}\log\frac{16}{\abs{f_t(a)-f_t(b)}}.
    \end{aligned}
\end{equation}

Put $K_t=\frac{1+t}{1-t}$ for $0\le t\le k$.
For $0<t<k$, $K_t>0$ and $K_t'=\frac{2}{(1-t)^2}=\frac{2K_t}{1-t^2}$. Put
\begin{equation}
    g(t) = K_t \log\frac{16}{\abs{f_t(a)-f_t(b)}}.
\end{equation}
The preceding differential inequality gives
\begin{equation}
    \label{eq:path-log-monotonicity}
    \begin{aligned}
        \frac{\dd}{\dd t}g(t) & =K_t'\log\frac{16}{\abs{f_t(a)-f_t(b)}} -K_t\frac{\dd}{\dd t}\log\abs{f_t(a)-f_t(b)} \\
                              & \ge \left(K_t'-\frac{2K_t}{1-t^2}\right) \log\frac{16}{\abs{f_t(a)-f_t(b)}}=0.
    \end{aligned}
\end{equation}
Therefore, $g$ is nondecreasing on $(0,k)$.
Since $k<1$, Theorem~\ref{thm:variation}, applied at $t=0$ and $t=k$, also gives continuity of $f_t(a)$ and $f_t(b)$ at both endpoints.
Since they remain distinct, $g$ is continuous on $[0,k]$, and hence it is nondecreasing there.
Using $f_0=\mathrm{id}$, $K_0=1$, and $K_k=\frac{1+k}{1-k}=K$, the inequality $g(0)\le g(k)$ gives
\begin{equation}
    \label{eq:path-endpoint-bound}
    K\log\frac{16}{\abs{f_k(a)-f_k(b)}}
    \ge\log\frac{16}{\abs{a-b}},
\end{equation}
which is equivalent to
\begin{equation}
    \abs{f_k(a) - f_k(b)} \le 16^{1-\frac{1}{K}} \abs{a-b}^{\frac{1}{K}}.
\end{equation}
Since $\mu_{f_k}=k\mu=\mu_f$, the composition formula \eqref{eq:qc-composition} shows that
$f\circ f_k^{-1}$ is a conformal automorphism of $\D$. It fixes $0$, so it is a rotation.
Thus \eqref{eq:main} follows.
Combining this upper bound with Proposition~\ref{prop:sharp} gives $M(K)=16^{1-\frac{1}{K}}$. \qed

\subsection{Sharpness}\label{subsec:sharpness}
It is well known that Mori's constant is at least $16^{1-\frac{1}{K}}$. For the reader's convenience, we include a proof here, based on the construction in \cite[Ch.~II, \S3.4, pp.~67--68]{LV}.
\begin{proposition}
    \label{prop:sharp}
    $M(K) \ge 16^{1-\frac{1}{K}}$ for every $K \ge 1$.
\end{proposition}

\begin{proof}
    Fix $K > 1$ and write $\D_+=\set{z\in\D:\im z>0}$. For $0<\theta<\pi$, put
    \begin{equation}
        q=\sin\frac{\theta}{2}, \quad q'=\cos\frac{\theta}{2},\quad \text{and}\quad m_\theta=\frac{\K(q')}{\K(q)}.
    \end{equation}
    Consider the quadrilateral
    $Q_\theta=(\D_+;0,1,\ee^{\mathrm{i}\theta},-1)$, with the vertices in the indicated order.
    It is well known that the Joukowski map $z \mapsto z+\frac{1}{z}$ sends $\D_+$ onto the lower half-plane. Therefore, the map
    \begin{equation}
        \label{eq:sharp-halfplane-map}
        T_\theta(z)=\frac{2-z-\frac{1}{z}}{4q^2}
    \end{equation}
    sends $\D_+$ onto the upper half-plane. Its boundary values at
    $1$, $\ee^{\mathrm{i}\theta}$, $-1$, and $0$ are $0$, $1$, $\frac{1}{q^2}$, and $\infty$, respectively.
    Using the Schwarz--Christoffel integral described in Section~\ref{subsec:elliptic}, define
    \begin{equation}
        \label{eq:sharp-rectangle-map}
        \Phi_\theta(z)=\frac{\K(q')}{\K(q)}
        +\frac{\mathrm{i}}{2\K(q)}\int_0^{T_\theta(z)}
        \frac{\dd u}{\sqrt{u(1-u)(1-q^2u)}}.
    \end{equation}
    Then $\Phi_\theta$ maps $\D_+$ conformally onto
    \begin{equation}
        \label{eq:sharp-rectangle}
        R_\theta=\set{x+\mathrm{i}y:0<x<m_\theta,\ 0<y<1}.
    \end{equation}
    It sends $0,1,\ee^{\mathrm{i}\theta},-1$ to
    $0,m_\theta,m_\theta+\mathrm{i},\mathrm{i}$, respectively.
    Therefore, the modulus of $Q_\theta$ is $m_\theta$, which is given by the classical formula
    \begin{equation}
        \label{eq:sharp-quadrilateral-modulus}
        m_\theta=\frac{2}{\pi}\mathcal{G}\left(\sin\frac{\theta}{2}\right).
    \end{equation}
    Choose $0<q_1<1$, set $q_2=\varphi_K(q_1)$, and put $\theta_j=2\arcsin q_j$ for $j=1,2$.
    The definition of $\varphi_K$ gives
    \begin{equation}
        \label{eq:sharp-modulus-relation}
        \mathcal{G}(q_1)=K\mathcal{G}(q_2),\qquad m_{\theta_1}=Km_{\theta_2}.
    \end{equation}
    Thus the affine map $A(x+\mathrm{i}y)=\frac{x}{K}+\mathrm{i}y$ is a $K$-quasiconformal map that takes $R_{\theta_1}$ onto $R_{\theta_2}$.
    Consequently,
    \begin{equation}
        \label{eq:sharp-semidisk-map}
        g(z)=\Phi_{\theta_2}^{-1}\circ A\circ\Phi_{\theta_1}(z) 
    \end{equation}
    is a $K$-quasiconformal homeomorphism from $Q_{\theta_1}$ to $Q_{\theta_2}$. It fixes $-1,0,1$ and sends
    $\ee^{\mathrm{i}\theta_1}$ to $\ee^{\mathrm{i}\theta_2}$. In particular, it maps the diameter $[-1,1]$ onto itself.
    Reflection in this diameter gives a $K$-quasiconformal map on the full disk:
    \begin{equation}
        \label{eq:sharp-disk-map}
        f(z)=
        \begin{cases}
            g(z),                  & \im z\ge0, \\
            \overline{g(\bar{z})}, & \im z<0.
        \end{cases}
    \end{equation}
    In particular, $f(\ee^{\pm\mathrm{i}\theta_1})=\ee^{\pm\mathrm{i}\theta_2}$.
    For the boundary points $a=\ee^{\mathrm{i}\theta_1}$ and $b=\ee^{-\mathrm{i}\theta_1}$, we have
    \begin{equation}
        \abs{a-b}=2\sin\theta_1=4q_1\sqrt{1-q_1^2}
    \end{equation}
    and
    \begin{equation}
        \abs{f(a)-f(b)} = 2 \sin\theta_2 =4q_2\sqrt{1-q_2^2}.
    \end{equation}
    The elliptic expansions
    \eqref{eq:K-series} and \eqref{eq:compl-series} give $\mathcal{G}(q)=\log\frac{4}{q}+o(1)$ as $q\to0$.
    Since $q_2\to0$ as $q_1\to0$, \eqref{eq:sharp-modulus-relation} implies
    \begin{equation}
        \label{eq:sharp-parameter-asymptotic}
        \log\frac{4}{q_2}=\frac{1}{K}\log\frac{4}{q_1}+o(1),\quad \text{and hence}\quad
        q_2=4^{1-\frac{1}{K}}q_1^{\frac{1}{K}}(1+o(1)).
    \end{equation}
    Therefore,
    \begin{equation}
        \label{eq:sharp-chord-limit}
        \frac{\abs{f(a)-f(b)}}{\abs{a-b}^{\frac{1}{K}}}
        =\frac{4q_2\sqrt{1-q_2^2}}{\left(4q_1\sqrt{1-q_1^2}\right)^{\frac{1}{K}}}
        \longrightarrow16^{1-\frac{1}{K}}\qquad(q_1\to0),
    \end{equation}
    and the conclusion follows.
\end{proof}

\section{Four linear bounds for the velocity}\label{sec:linear}

Throughout this section, $a,b \in \D $ satisfy $a \ne b$ and $ \abs {a} \ge \abs {b}$.
We assume $b \ne 0$.
Also, $ \norm { \mu }_ \infty \le 1$, $f_t$ is the reflected family of $ \mu $, and $v$ is its initial velocity field
(Proposition~\ref{prop:velocity}). Write
\begin{equation}
    \label{eq:pair}
    \begin{gathered}
        s= \frac { \abs {a-b}}{ \abs {a}+ \abs {b}} \in (0,1], \qquad
        \eta = \frac { \abs {a}- \abs {b}}{ \abs {a}+ \abs {b}} \in [0,1), \\
        q_a= \frac { \abs {a}^2}{1- \abs {a}^2}, \qquad q_b= \frac { \abs {b}^2}{1- \abs {b}^2}.
    \end{gathered}
\end{equation}
We call $s$ the \emph{relative chord} and write $\delta=\delta(a,b)$ for the pseudo-hyperbolic
distance defined in Section~\ref{subsec:distortion}. The velocity enters through the real-linear functionals
\begin{equation}
    \label{eq:functionals}
    \begin{aligned}
        u_a & = \re \frac {v(a)}{a}, \qquad u_b= \re \frac {v(b)}{b},                           \\
        \Lc & = \re \frac {v(a)-v(b)}{a-b}, \qquad W= \frac {v(a)}{a}- \frac {v(b)}{b} \in \C .
    \end{aligned}
\end{equation}
By Lemma~\ref{lem:log-modulus-derivative}, $u_a= \frac { \dd }{ \dd t} \log \abs {f_t(a)}$ and $ \Lc = \frac { \dd }{ \dd t} \log \abs {f_t(a)-f_t(b)}$ at $t=0$. For $0<s \le 1$ put
\begin{equation}
    \label{eq:Aang}
    p_s= \sin \Bigl ( \frac {1}{2} \arcsin s \Bigr ), \qquad \A (s)= \sqrt {1-s^2} \; B(p_s),
    \qquad \A (1)=0,
\end{equation}
and let $ \J (a,b)= \int _ \C \frac { \dA (w)}{ \abs {w} \, \abs {w-a} \, \abs {w-b}}$ (studied in Section~\ref{sec:threepole}).

\begin{proposition}[Four linear bounds]
    \label{prop:four}
    For every $ \mu $ with $ \norm { \mu }_ \infty \le 1$:
    \begin{equation}
        \label{eq:R}
        \makebox[0.7\linewidth][l]{$\displaystyle
                \text{\makebox[2.5em][l]{\textup{(R)}}}\abs{u_a}\le2m(\abs{a}),\qquad\abs{u_b}\le2m(\abs{b});
            $}
    \end{equation}
    \begin{equation}
        \label{eq:S}
        \makebox[0.7\linewidth][l]{$\displaystyle
                \text{\makebox[2.5em][l]{\textup{(S)}}}\Lc+q_au_a+q_bu_b\le B(\delta);
            $}
    \end{equation}
    \begin{equation}
        \label{eq:A}
        \makebox[0.7\linewidth][l]{$\displaystyle
                \text{\makebox[2.5em][l]{\textup{(A)}}}\Lc-\frac{\abs{a}u_a+\abs{b}u_b}{\abs{a}+\abs{b}}\le\A(s);
            $}
    \end{equation}
    \begin{equation}
        \label{eq:W}
        \makebox[0.7\linewidth][l]{$\displaystyle
                \text{\makebox[2.5em][l]{\textup{(W)}}}\abs{W}\le\frac{\abs{a-b}}{\pi}\,\J(a,b).
            $}
    \end{equation}
\end{proposition}

\begin{proof}
    For $0 \le t<1$ the map $f_t$ is $K_t$-quasiconformal with $K_t= \frac {1+t}{1-t}$, so $K_t'(0)=2$. Each of (R),
    (S), (A) comes from a distortion inequality that is an equality at $t=0$: the right derivatives at $t=0$
    of the two sides are therefore ordered.

    (R) By Theorem~\ref{thm:schwarz}, $ \abs {f_t(a)} \le \varphi _{K_t}( \abs {a})$. Lemma~\ref{lem:elliptic}(a) gives
    $u_a \le 2m( \abs {a})$. Replacing $ \mu $ by $- \mu $ (which replaces $f_t$ by $f_{-t}$ and $v$ by $-v$) gives
    $-u_a \le 2m( \abs {a})$. The same argument applies to $b$.
    This radial bound also follows from Semenov \cite[(4)--(6), p.~504]{Sem91},
    whose normalisation $\abs{u_{\bar{z}}}\le\frac{1}{2}$ corresponds to $u=\frac{v}{2}$.

    (S) Let $S_t= \frac{\abs {f_t(a)-f_t(b)}}{\sqrt {(1- \abs {f_t(a)}^2)(1- \abs {f_t(b)}^2)}}$. Then $S_t= \frac{\delta _t}{\sqrt {1- \delta _t^2}}$,
    where $ \delta _t= \delta (f_t(a),f_t(b))$, and $S_t$ is an increasing function of $ \delta _t$. Since
    $ \frac { \dd }{ \dd t} \log (1- \abs {f_t(a)}^2)=-2q_au_a$ at $t=0$,
    \begin{equation}
        \frac { \dd }{ \dd t} \log S_t \Big |_{t=0}= \Lc +q_au_a+q_bu_b.
    \end{equation}
    On the other hand $ \delta _t \le \varphi _{K_t}( \delta )$ by Theorem~\ref{thm:schwarz}, and by
    Lemma~\ref{lem:elliptic}(a)
    \begin{equation}
        \frac { \dd }{ \dd t} \log \frac { \varphi _{K_t}( \delta )}{ \sqrt {1- \varphi _{K_t}( \delta )^2}} \Big |_{t=0}
        = \frac {2m( \delta )}{1- \delta ^2}=B( \delta ).
    \end{equation}

    (A) Apply Theorem~\ref{thm:angle} to $G=F_t$ (which fixes $ \infty $) and the triple $(0,a,b)$. Here
    $ \sin \frac {\Theta} {2}=s$, so $ \sin \frac {\Theta} {4}=p_s$. With
    \begin{equation}
        T_t= \frac { \abs {f_t(a)-f_t(b)}}{ \abs {f_t(a)}+ \abs {f_t(b)}}, \qquad \text {so that} \qquad
        \frac { \dd }{ \dd t} \log T_t \Big |_{t=0}= \Lc - \frac { \abs {a}u_a+ \abs {b}u_b}{ \abs {a}+ \abs {b}},
    \end{equation}
    the theorem reads $p_{T_t} \le \varphi _{K_t}(p_s)$. If $s=1$, then $T_t \le 1=T_0$ by the triangle
    inequality, which gives (A) with $ \A (1)=0$. If $0<s<1$, write $p=p_s$. Then
    \begin{equation}
        s=2p \sqrt {1-p^2}, \qquad \frac { \dd \log p_s}{ \dd \log s}= \frac {1-p^2}{ \sqrt {1-s^2}},
    \end{equation}
    so differentiation at $t=0$ gives
    \begin{equation}
        \frac {1-p^2}{ \sqrt {1-s^2}} \Bigl ( \Lc - \frac { \abs {a}u_a+ \abs {b}u_b}{ \abs {a}+ \abs {b}} \Bigr ) \le 2m(p)=(1-p^2)B(p),
    \end{equation}
    which is (A).

    (W) The kernel $ \mathcal {R}(w,z)= \frac {z(z-1)}{w(w-1)(w-z)}$ of \eqref{eq:V} satisfies
    \begin{equation}
        \frac { \mathcal {R}(w,a)}{a}- \frac { \mathcal {R}(w,b)}{b}= \frac {a-b}{w(w-a)(w-b)}.
    \end{equation}
    Since $v=V[ \mu ^{ \mathrm {e}}]$ on $ \D $, this gives
    \begin{equation}
        W=- \frac {a-b} {\pi} \int _ \C \frac { \mu ^{ \mathrm {e}}(w)}{w(w-a)(w-b)} \dA (w),
    \end{equation}
    and $ \abs { \mu ^{ \mathrm {e}}} \le 1$ on $ \C $ gives (W).
\end{proof}

\begin{remark}
    Each of (R)--(W) is a separate inequality valid for every $ \mu $. The rest of the proof only forms
    nonnegative combinations of them; it never needs the bounds to be attained by the same $ \mu $.
\end{remark}

\section{A three-pole integral}\label{sec:threepole}

For distinct nonzero $a,b \in \C $,
\begin{equation}
    \label{eq:J}
    \J (a,b)= \int _ \C \frac { \dA (w)}{ \abs {w} \, \abs {w-a} \, \abs {w-b}}
\end{equation}
is finite: the three singularities are integrable and the integrand is $O( \abs {w}^{-3})$ at $ \infty $. By
congruence invariance it depends only on the three lengths $ \abs {a}$, $ \abs {b}$, $ \abs {a-b}$.
When $ \J $ is given three positive real arguments, these denote the three side lengths of the corresponding triangle (degenerate triangles allowed).
Therefore, we also write
\begin{equation}
    \J ( \abs {a}, \abs {b}, \abs {a-b})= \J (a,b).
\end{equation}
The integral $ \J $ is homogeneous of degree $-1$. For each $t >0$,
\begin{equation}
    \J (ta,tb)=t^{-1} \J (a,b) = t^{-1} \J ( \abs {a}, \abs {b}, \abs {a-b}).
\end{equation}
Since $ \J (a,b)=\J (b,a)$, we use the same ordering convention $ \abs {a} \ge \abs {b}>0$
in the proofs in this section. With $s$ and $ \eta $ as in \eqref{eq:pair}, put
\begin{equation}
    \J _ \eta (s)= \J (1+ \eta , \, 1- \eta , \, 2s), \qquad \Jz (s)= \J (1,1,2s),
\end{equation}
and then we have
\begin{equation}
    \label{eq:jabjeta}
    \begin{aligned}
        \J (a,b)= & \frac{2}{\abs {a}+ \abs {b}} \J\left(\frac{2\abs{a}}{\abs{a}+\abs{b}}, \frac{2\abs{b}}{\abs{a}+\abs{b}}, \frac{2\abs{a-b}}{\abs{a}+\abs{b}}\right) \\
        =         & \frac {2 \, \J _ \eta (s)}{ \abs {a}+ \abs {b}} = \frac { \sqrt {1- \eta ^2}}{ \sqrt { \abs {ab}}} \J _ \eta (s).
    \end{aligned}
\end{equation}

The integrand defining $\J(a,b)$ has singularities at $0,a,b$, making the integral difficult to estimate directly. The Feynman-parameter identity (see Lemma~\ref{lem:feynman}) rewrites it as an integral over a triangle.
\begin{proposition}
    \label{prop:simplex-representation}
    Let $a,b \in \C $ be distinct and nonzero. Set
    \begin{equation}
        \Sigma = \set {( \xi , \zeta ) \in \R ^2: \xi >0, \  \zeta >0, \  \xi + \zeta <1}
        \quad \text{and} \quad  \ell =1- \xi - \zeta .
    \end{equation}
    Then we have
    \begin{equation}
        \label{eq:simplex}
        \J (a,b)= \int _ \Sigma \frac { \dd \xi \dd \zeta }{ \sqrt { \xi \zeta \ell \, ( \xi \zeta \abs {a}^2+ \xi \ell \abs {b}^2+ \zeta \ell \abs {a-b}^2)}}.
    \end{equation}
\end{proposition}

\begin{proof}
    For $w \notin \set {0,a,b}$, substitute $A= \abs {w}^2$, $B= \abs {w-a}^2$, $C= \abs {w-b}^2$
    into \eqref{eq:feynman}. Since the resulting integrand is nonnegative,
    Tonelli's theorem allows us to exchange the order of integration and obtain
    \begin{equation}
        \J (a,b)= \frac {1}{2 \pi } \int _ \Sigma \frac {1}{ \sqrt { \xi \zeta \ell }}
        \left ( \int _ \C \frac { \dA (w)}{( \xi \abs {w}^2+ \zeta \abs {w-a}^2+ \ell \abs {w-b}^2)^{\frac{3}{2}}} \right ) \dd \xi \dd \zeta .
    \end{equation}
    Set $c= \zeta a+ \ell b$. Since $ \xi + \zeta + \ell =1$, completing the square gives
    \begin{equation}
        \begin{aligned}
            \xi \abs {w}^2+ \zeta \abs {w-a}^2+ \ell \abs {w-b}^2                                               = & \abs{w}^2 - 2 \re \left(\left(\overline{\zeta a + \ell b} \right)w\right) + \zeta \abs{a}^2 + \ell \abs{b}^2 \\
            =                                                                                                     & \abs {w-c}^2+ \zeta \abs {a}^2+ \ell \abs {b}^2- \abs {c}^2                                                  \\
            =                                                                                                     & \abs {w-c}^2+ \Delta.
        \end{aligned}
    \end{equation}
    Expanding $c = \zeta a + \ell b$ and using $2\re (\bar{a} b) = \abs{a}^2 + \abs{b}^2 - \abs{a-b}^2$ gives
    \begin{equation} \label{eq:Delta}
        \begin{aligned}
            \Delta & = \zeta (1- \zeta ) \abs {a}^2+ \ell (1- \ell ) \abs {b}^2-2 \zeta \ell \re ( \ol {a}b)                        \\
                   & = \zeta (1- \zeta ) \abs {a}^2+ \ell (1- \ell ) \abs {b}^2- \zeta \ell ( \abs {a}^2+ \abs {b}^2- \abs {a-b}^2) \\
                   & = \xi \zeta \abs {a}^2+ \xi \ell \abs {b}^2+ \zeta \ell \abs {a-b}^2 > 0.
        \end{aligned}
    \end{equation}
    Applying Lemma~\ref{lem:plane-radial-integral} with centre $c$ gives
    \begin{equation}
        \J (a,b)= \frac {1}{2 \pi } \int _ \Sigma \frac {1}{ \sqrt { \xi \zeta \ell }}
        \left ( \int _ \C \frac { \dA (w)}{(\abs{w-c}^2 + \Delta)^{\frac{3}{2}}} \right ) \dd \xi \dd \zeta =
        \int _ \Sigma \frac {1}{ \sqrt { \xi \zeta \ell }\sqrt{\Delta}} \dd \xi \dd \zeta .
    \end{equation}
\end{proof}

\begin{lemma}
    \label{lem:simplex}
    For distinct nonzero $a,b \in \C $,
    \begin{equation} \label{eq:simplex-comparison}
        \sqrt { \abs {ab}} \; \J (a,b) \le \Jz  (s), \quad \text{with} \quad s =  \frac { \abs {a-b}}{ \abs {a}+ \abs {b}},
    \end{equation}
    and the strict inequality holds if $ \abs {a} \ne \abs {b}$. Equivalently,
    \begin{equation} \label{eq:simplex2}
        \sqrt {1- \eta ^2} \, \J _ \eta (s) \le \Jz (s), \quad \text{with} \quad \eta =  \frac { \abs {a}- \abs {b}}{ \abs {a}+ \abs {b}}.
    \end{equation}
\end{lemma}

\begin{proof}
    By symmetry, we may assume $ \abs {a} \ge \abs {b}>0$. Then $0 \le \eta <1$,
    and the triangle inequality gives $0<s \le 1$ and $\eta \le s$.
    By \eqref{eq:jabjeta}, $ \sqrt { \abs {ab}} \, \J (a,b)= \sqrt {1- \eta ^2} \, \J _ \eta (s)$.
    It therefore suffices to compare $ \sqrt {1- \eta ^2} \, \J _ \eta (s)$ with $ \Jz (s)$.

    The substitution
    \begin{equation}
        \xi =1-u, \qquad \zeta = \frac {u(1+v)}{2}, \qquad \ell = \frac {u(1-v)}{2}
        \qquad (0<u<1, \  -1<v<1)
    \end{equation}
    maps the rectangle onto $ \Sigma $ and has absolute Jacobian $\frac{u}{2}$.
    Put $r=\frac{\abs {a}+ \abs {b}}{2}>0$, so that $ \abs {a}=r(1+ \eta )$, $ \abs {b}=r(1- \eta )$, and $ \abs {a-b}=2rs$.
    Substituting these expressions into $\Delta$, which is given in \eqref{eq:Delta}, gives
    \begin{equation}
        \label{eq:simplex-rectangle}
        \begin{aligned}
            \Delta & = (1-u)\frac{u(1+v)}{2}r^2(1+\eta)^2 +(1-u)\frac{u(1-v)}{2}r^2(1-\eta)^2
            +\frac{u^2(1-v^2)}{4}(2rs)^2                                                                        \\
                   & = \frac {r^2u(1-u)}{2} \Bigl ((1+v)(1+ \eta )^2+(1-v)(1- \eta )^2 \Bigr )+r^2u^2(1-v^2)s^2 \\
                   & = \frac{r^2u(1-u)}{2}\Bigl(2(1+\eta^2)+4\eta v\Bigr)+r^2u^2(1-v^2)s^2                      \\
                   & =r^2 \Bigl (u(1-u)(1+ \eta ^2+2 \eta v)+u^2(1-v^2)s^2 \Bigr ).
        \end{aligned}
    \end{equation}
    Consequently,
    \begin{equation}
        \label{eq:Jeta-rectangle}
        \begin{aligned}
            \J _ \eta (s) & =r \, \J (a,b)                                                                                                                \\
                          & = \int _0^1 \! \! \int _{-1}^1 \frac { \dd v \dd u}{ \sqrt {(1-u)(1-v^2)} \sqrt {u(1-u)(1+ \eta ^2+2 \eta v)+u^2(1-v^2)s^2}}.
        \end{aligned}
    \end{equation}
    We now make a further change of variables, by setting $v=\frac{q-\eta}{1-\eta q}$.
    Because $0 \le \eta <1$, this is an increasing bijection
    from $-1<q<1$ onto $-1<v<1$. Direct calculation gives
    \begin{equation}
        \begin{aligned}
             & \frac { \dd v}{ \dd q}  = \frac {1- \eta ^2}{(1- \eta q)^2}, \qquad 1-v^2= \frac {(1- \eta ^2)(1-q^2)}{(1- \eta q)^2}, \\
             & 1+ \eta ^2+2 \eta v     = \frac {(1- \eta ^2)(1+ \eta q)}{1- \eta q}.
        \end{aligned}
    \end{equation}
    Substituting these expressions into \eqref{eq:Jeta-rectangle} gives
    \begin{equation}
        \label{eq:simplex-normalised}
        \sqrt {1- \eta ^2} \, \J _ \eta (s)= \int _0^1 \! \! \int _{-1}^1
        \frac { \sqrt {1- \eta ^2} \dd q \dd u}{ \sqrt {1-u} \sqrt {1-q^2} \sqrt { \Theta _ \eta }},
    \end{equation}
    where
    \begin{equation}
        \label{eq:delta-eta}
        \Theta _ \eta =u(1-u)(1- \eta ^2q^2)+u^2(1-q^2)s^2>0.
    \end{equation}
    Setting $\eta=0$ in \eqref{eq:delta-eta} gives $\Theta_0=u(1-u)+u^2(1-q^2)s^2$. Therefore,
    \begin{equation}
        \begin{aligned}
              & \Theta_\eta-(1-\eta^2)\Theta_0                                                    \\
            = & u(1-u)(1-\eta^2q^2)+u^2(1-q^2)s^2 -(1-\eta^2)\Bigl(u(1-u)+u^2(1-q^2)s^2\Bigr)     \\
            = & u(1-u)\Bigl((1-\eta^2q^2)-(1-\eta^2)\Bigr) +u^2(1-q^2)s^2\Bigl(1-(1-\eta^2)\Bigr) \\
            = & \eta^2(1-q^2)u(1-u)+\eta^2u^2(1-q^2)s^2 \ge 0,
        \end{aligned}
    \end{equation}
    and it is strict whenever $\eta\ne 0$. Hence,
    \begin{equation}
        \frac { \sqrt {1- \eta ^2}}{ \sqrt { \Theta _ \eta }} \le \frac {1}{ \sqrt { \Theta _0}}.
    \end{equation}
    Therefore, inequality \eqref{eq:simplex2} holds, and the strict inequality holds whenever $ \eta \ne 0$.
\end{proof}

\begin{proposition}[The series expansion for $ \Jz $]
    \label{prop:J0}
    For $0<s \le 1$,
    \begin{equation}
        \label{eq:J0-elliptic}
        \Jz (s)=4 \int _0^{ \frac{\pi}{2}} \K \bigl ( \sqrt {1-s^2 \sin ^2 \theta } \bigr ) \dd \theta .
    \end{equation}
    For $0<s<1$,
    \begin{equation}
        \label{eq:J0-series}
        \frac { \Jz (s)}{ \pi }=2 \sum _{n \ge 0}b_n^3s^{2n} \Bigl ( \log \frac {8}{s}-3d_n \Bigr ).
    \end{equation}
    Consequently,
    \begin{equation}
        \label{eq:J0-bounds}
        \frac { \Jz (s)}{ \pi } \le 3+(2 \log \frac {8}{s}-3)(1-s^2)^{-\frac{1}{8}}.
    \end{equation}
    In particular,
    \begin{equation}
        \label{eq:J0-coarse}
        \frac { \Jz (s)}{ \pi } \le 2 \log \frac {8}{s} \, (1-s^2)^{-\frac{1}{8}}.
    \end{equation}
\end{proposition}

\begin{proof}
    At $ \eta =0$, we have $ \Theta _0=u \bigl (1-u+u(1-q^2)s^2 \bigr )$, so
    \eqref{eq:simplex-normalised} becomes
    \begin{equation}
        \label{eq:J0-iterated}
        \Jz (s)= \int _0^1 \! \! \int _{-1}^{1}
        \frac { \dd q \dd u}{\sqrt {1-u} \sqrt {1-q^2} \sqrt {u(1-u+u(1-q^2)s^2)}}.
    \end{equation}
    Set
    \begin{equation}
        u= \sin ^2 \varphi , \qquad q= \cos \theta ,
        \qquad 0< \varphi < \frac { \pi }{2}, \quad 0< \theta < \pi .
    \end{equation}
    This is a smooth bijection onto $(0,1) \times (-1,1)$, with absolute Jacobian
    \begin{equation}
        \abs { \frac { \partial (u,q)}{ \partial ( \varphi , \theta )}}
        = \abs {-2 \sin \varphi \cos \varphi \sin \theta }
        =2 \sin \varphi \cos \varphi \sin \theta .
    \end{equation}
    Substituting these expressions into \eqref{eq:J0-iterated} gives
    \begin{equation}
        \begin{aligned}
            \Jz (s) & =  \int_0^ \pi \int _0^{ \frac{\pi}{2}} \frac{2\sin\varphi \cos\varphi \sin\theta}{\sin\varphi \cos\varphi \cdot \sin\theta \cdot \sqrt{1-(1-s^2 \sin^2\theta)\sin^2\varphi}} \dd\varphi\dd\theta \\
                    & =2 \int _0^ \pi \int _0^{ \frac{\pi}{2}}
            \frac { \dd \varphi \dd \theta }{ \sqrt {1-(1-s^2 \sin ^2 \theta ) \sin ^2 \varphi }}                                                                                                                       \\
                    & =2 \int _0^ \pi \K \bigl ( \sqrt {1-s^2 \sin ^2 \theta } \bigr ) \dd \theta  =4 \int _0^{ \frac{\pi}{2}} \K \bigl ( \sqrt {1-s^2 \sin ^2 \theta } \bigr ) \dd \theta .
        \end{aligned}
    \end{equation}
    For $0<s<1$ and $0<\theta<\frac{\pi}{2}$, we have $0<s\sin\theta<1$ and, since $d_n<\log2$,
    $\log(\frac{4}{s\sin\theta})-2d_n>0$ for every $n\ge0$.
    Thus the series obtained by substituting $s\sin\theta$ for $q$ in \eqref{eq:compl-series}
    has nonnegative terms. Tonelli's theorem permits termwise integration, and the two integral
    identities in \eqref{eq:sine-moments} then give \eqref{eq:J0-series}:
    \begin{equation}
        \begin{aligned}
            \J_0(s) & = 4 \sum _{n \ge 0} \int _0^{ \frac{\pi}{2}}b_n^2s^{2n} \sin ^{2n} \theta \Bigl ( \log \frac {4}{s \sin \theta }-2d_n \Bigr ) \dd \theta                                                             \\
                    & = 4 \sum_{n \ge 0} b^2_n s^{2n} \left(\left(\log \frac{4}{s}- 2 d_n\right) \int^{\frac{\pi}{2}}_0 \sin^{2n} \theta \, d\theta - \int^{\frac{\pi}{2}}_0 \sin^{2n} \theta \log(\sin \theta) \, d\theta
            \right)                                                                                                                                                                                                        \\
                    & =2 \pi \sum _{n \ge 0}b_n^3s^{2n} \Bigl ( \log \frac {4}{s}-2d_n+ \log 2-d_n \Bigr )                                                                                                                 \\
                    & =2 \pi \sum _{n \ge 0}b_n^3s^{2n} \Bigl ( \log \frac {8}{s}-3d_n \Bigr ).
        \end{aligned}
    \end{equation}

    \par Let $F(s^2)= \sum _{n \ge 0}b_n^3s^{2n}$. We claim that
    \begin{equation} \label{eq:Fs2}
        b_n^3 \le {(\frac{1}{8})_n}/{n!}.
    \end{equation}
    Both sides are $1$ at $n=0$, and the ratio of
    consecutive terms is $\left(\frac{b_{n+1}}{b_n}\right)^3 = \bigl ( \frac {n+\frac{1}{2}}{n+1} \bigr )^3$ on the left and $ \frac {n+\frac{1}{8}}{n+1}$ on the right. Direct calculation shows that $\bigl ( \frac {n+\frac{1}{2}}{n+1} \bigr )^3 \le  \frac {n+\frac{1}{8}}{n+1}$. Therefore, the claim follows. Combining this inequality with \eqref{eq:binomial} gives
    \begin{equation} \label{eq:Fs2-bound}
        F(s^2)= \sum _{n \ge 0}b_n^3s^{2n} \le \sum _{n \ge 0} \frac {(\frac{1}{8})_n}{n!}s^{2n}=(1-s^2)^{-\frac{1}{8}}.
    \end{equation}
    In addition, since $d_n \ge d_1= \frac {1}{2}$ for $n \ge 1$, there holds
    \begin{equation}
        \sum_{n \ge 0} d_n b_n^3 s^{2n} \ge \frac{1}{2} \sum_{n \ge 0}b_n^3 s^{2n} - \frac{1}{2} b_0^3 = \frac{1}{2}F(s^2) - \frac{1}{2}.
    \end{equation}
    Therefore,
    \begin{equation} \label{eq:J0-bound2}
        \frac{\Jz (s)}{\pi}  \le  2 F(s^2) \log \frac{8}{s} - 3(F(s^2) - 1) = 3 + (2 \log \frac{8}{s} - 3) F(s^2).
    \end{equation}
    Since $2\log(\frac{8}{s})-3\ge2\log8-3>0$,
    substituting \eqref{eq:Fs2-bound} into \eqref{eq:J0-bound2} yields \eqref{eq:J0-bounds}.
    Finally, \eqref{eq:J0-coarse} follows from $3 \bigl ((1-s^2)^{-\frac{1}{8}}-1 \bigr ) \ge 0$.
\end{proof}

\section{Two bounds for \texorpdfstring{$ \Lc $}{L}}\label{sec:joint}

Throughout this section, we continue to assume $a,b \in \D $, $a \ne b$, and $ \abs {a} \ge \abs {b}>0$.
We retain $s$ and $ \eta $ from \eqref{eq:pair}. In addition, we use
\begin{equation}
    \label{eq:joint-notation}
    \begin{gathered}
        \varrho = \frac { \abs {b}}{ \abs {a}} \in (0,1], \qquad D=(1- \abs {a}^2)(1- \abs {b}^2), \\
        M= \abs {a}^2+ \abs {b}^2-2 \abs {ab}^2, \qquad \omega = \frac {D}{1- \abs {ab}^2} \in (0,1), \\
        \Hc =2 \log \frac {16}{ \abs {a-b}}=2 \log \frac {8}{s}+2 \log \frac {2}{ \abs {a}+ \abs {b}}.
    \end{gathered}
\end{equation}
Then $ \eta =\frac{1- \varrho }{1+ \varrho } \in [0,1)$, and the triangle inequality gives $ \eta ^2 \le s^2 \le 1$. A direct computation gives
\begin{equation}
    \label{eq:Z}
    \frac {1}{ \omega }=1+q_a+q_b= \frac {1- \abs {ab}^2}{D}, \qquad 1- \omega = \frac {M}{1- \abs {ab}^2}.
\end{equation}
Our goal is $ \Lc \le \Hc $ for every direction $ \mu $ with $ \norm { \mu }_ \infty \le 1$, as asserted in
Theorem~\ref{thm:kernel}.

\subsection{The first bound}
The functionals $u_a,u_b$ are linked through $W$:
\begin{equation}
    \label{eq:L-identities}
    \Lc =u_a+ \re \frac {bW}{a-b}=u_b+ \re \frac {aW}{a-b}.
\end{equation}
Indeed, $ \frac {v(a)-v(b)}{a-b}= \frac {v(a)}{a}+ \frac {b}{a-b}W= \frac {v(b)}{b}+ \frac {a}{a-b}W$. Multiplying the two
identities in \eqref{eq:L-identities} by $q_a$ and $q_b$ and adding $ \Lc $ gives the exact identity
\begin{equation}
    \label{eq:ZL}
    \frac { \Lc }{ \omega }= \bigl ( \Lc +q_au_a+q_bu_b \bigr )+ \re \Bigl ( \frac {q_ba+q_ab}{a-b} \, W \Bigr ).
\end{equation}

\begin{proposition}[First bound]
    \label{prop:joint1}
    With
    \begin{equation}
        \label{eq:h}
        h= \frac { \omega \, \abs {q_ba+q_ab}}{ \sqrt { \abs {ab}}},
    \end{equation}
    we have
    \begin{equation}
        \label{eq:joint1}
        \Lc \le \omega \Bigl (B( \delta )+ \abs {q_ba+q_ab} \, \frac{\J (a,b)}{\pi} \Bigr ) \le \omega B( \delta )+h \, \frac { \Jz (s)}{ \pi }.
    \end{equation}
\end{proposition}

\begin{proof}
    By (S) in Proposition~\ref{prop:four}, the first term on the right side of \eqref{eq:ZL} satisfies
    $\Lc+q_au_a+q_bu_b\le B(\delta)$.
    For the second term, use (W) to obtain
    \begin{equation}
        \label{eq:joint1-W-estimate}
        \begin{aligned}
            \re\Bigl(\frac{q_ba+q_ab}{a-b}W\Bigr)
            \le \frac{\abs{q_ba+q_ab}}{\abs{a-b}}\abs{W}
            \le \frac{\abs{q_ba+q_ab}}{\pi}\J(a,b).
        \end{aligned}
    \end{equation}
    Substituting these two estimates into \eqref{eq:ZL} gives
    \begin{equation}
        \label{eq:joint1-before-comparison}
        \frac{\Lc}{\omega}
        \le B(\delta)+\frac{\abs{q_ba+q_ab}}{\pi}\J(a,b).
    \end{equation}
    Multiplying by $\omega>0$ proves the first inequality in \eqref{eq:joint1}.
    Together with Lemma~\ref{lem:simplex} and the definition of $h$, we obtain
    \begin{equation}
        \label{eq:joint1-integral-comparison}
        \omega\frac{\abs{q_ba+q_ab}}{\pi}\J(a,b)
        \le \frac{\omega\abs{q_ba+q_ab}}{\sqrt{\abs{ab}}}\frac{\Jz(s)}{\pi}
        =h\frac{\Jz(s)}{\pi}.
    \end{equation}
    Combining this with the first inequality proves the second inequality in \eqref{eq:joint1}.
\end{proof}

\subsection{The coefficient $h$}
Let $ \alpha = \varrho + \varrho ^{-1} \ge 2$ and define $ \beta \ge 0$ and $ \gamma $ by
\begin{equation}
    \label{eq:beta}
    \beta ^2= \frac { \abs {ab}(1+ \abs {ab})^2( \abs {a}- \abs {b})^2}{M^2}, \qquad \gamma = \frac {( \alpha +2)D}{( \alpha -2 \abs {ab})^2}.
\end{equation}

\begin{lemma}
    \label{lem:coefficient}
    \begin{equation}
        \label{eq:beta-alpha}
        \beta ^2= \frac {( \alpha -2)(1+ \abs {ab})^2}{( \alpha -2 \abs {ab})^2},
    \end{equation}
    \begin{equation}
        \label{eq:beta-gamma-bound}
        \beta ^2+ \gamma =1- \frac {( \alpha -2)(1+ \abs {ab})}{ \alpha -2 \abs {ab}} \le 1,
    \end{equation}
    \begin{equation}
        \label{eq:h-identity}
        \Bigl ( \frac {h}{1- \omega } \Bigr )^2= \beta ^2+ \gamma \, (1-s^2),
    \end{equation}
    and
    \begin{equation}
        \label{eq:h-shape}
        \frac {h}{1- \omega } \le \max \set { \sqrt {1-s^2}, \  \sqrt { \varrho }}.
    \end{equation}
\end{lemma}

\begin{proof}
    Direct computation shows that
    \begin{equation} \label{eq:joint-notation2}
        \begin{aligned}
            ( \abs {a}+ \abs {b})^2= \abs {ab}( \alpha +2), \qquad & ( \abs {a}- \abs {b})^2= \abs {ab}( \alpha -2), \quad M= \abs {ab}( \alpha -2 \abs {ab}), \\
            D=1+ \abs {ab}^2- \abs {ab} \alpha, \quad              & \abs{a-b}^2 =(\abs{a}+ \abs{b})^2 s^2 = s^2\abs{ab}( \alpha +2).
        \end{aligned}
    \end{equation}
    Substituting these equations into the definition of $\beta^2$ gives \eqref{eq:beta-alpha}.
    Moreover, \eqref{eq:beta-gamma-bound} follows from the following identity
    \begin{equation}
        \label{eq:beta-gamma-numerator}
        \begin{aligned}
                & ( \alpha -2)(1+ \abs {ab})^2+( \alpha +2)D                                            \\
            ={} & ( \alpha -2)(1+2 \abs {ab}+ \abs {ab}^2)
            +( \alpha +2)(1+ \abs {ab}^2- \alpha \abs {ab})                                             \\
            ={} & 2 \alpha - \alpha ^2 \abs{ab} - 4\abs {ab}+2 \alpha \abs {ab}^2                       \\
            ={} & ( \alpha -2 \abs {ab})(2- \alpha \abs {ab})                                           \\
            ={} & ( \alpha -2 \abs {ab}) \bigl(( \alpha -2 \abs {ab})-( \alpha -2)(1+ \abs {ab})\bigr).
        \end{aligned}
    \end{equation}
    Now we prove \eqref{eq:h-identity}. By \eqref{eq:Z},
    \begin{equation} \label{eq:h-omega}
        \frac {h}{1- \omega }= \frac {D \, \abs {q_ba+q_ab}}{ \sqrt { \abs {ab}} \, M}.
    \end{equation}
    Expanding with $2 \re (a \bar {b})= \abs {a}^2+ \abs {b}^2- \abs {a-b}^2$, we have
    \begin{equation}
        \label{eq:h-numerator-expand}
        \begin{aligned}
            \frac {D^2\abs {q_ba+q_ab}^2}{ \abs {ab}^2}
             & = \abs {(1- \abs {a}^2)\frac{a \abs{b}}{\abs{a}}+(1- \abs {b}^2) \frac{\abs{a}b}{\abs{b}}}^2 \\
             & = \abs {b}^2(1- \abs {a}^2)^2+ \abs {a}^2(1- \abs {b}^2)^2+2D \re (a \bar {b})               \\
             & = \abs {b}^2(1- \abs {a}^2)^2+ \abs {a}^2(1- \abs {b}^2)^2
            +D \bigl ( \abs {a}^2+ \abs {b}^2- \abs {a-b}^2 \bigr ).
        \end{aligned}
    \end{equation}
    To combine the first two terms, complete the square and use
    $ \abs {b}(1- \abs {a}^2)+ \abs {a}(1- \abs {b}^2)=( \abs {a}+ \abs {b})(1- \abs {ab})$:
    \begin{equation}
        \label{eq:h-numerator-square}
        \begin{aligned}
            \abs {b}^2(1- \abs {a}^2)^2+ \abs {a}^2(1- \abs {b}^2)^2 & = \bigl ( \abs {b}(1- \abs {a}^2)+ \abs {a}(1- \abs {b}^2) \bigr )^2-2 \abs {ab}D \\
                                                                     & =( \abs {a}+ \abs {b})^2(1- \abs {ab})^2-2 \abs {ab}D.
        \end{aligned}
    \end{equation}
    Substituting \eqref{eq:h-numerator-square} into \eqref{eq:h-numerator-expand} and collecting the terms containing $D$, we obtain
    \begin{equation}
        \label{eq:h-numerator}
        \begin{aligned}
            \frac {D^2\abs{q_ba+q_ab}^2}{ \abs {ab}^2} & =( \abs {a}+ \abs {b})^2(1- \abs {ab})^2  +D \bigl ( \abs {a}^2+ \abs {b}^2-2 \abs {ab}- \abs {a-b}^2 \bigr ) \\
                                                       & =( \abs {a}+ \abs {b})^2(1- \abs {ab})^2 -D \bigl ( \abs {a-b}^2-( \abs {a}- \abs {b})^2 \bigr )              \\
                                                       & = \abs{ab} (\alpha+2) (1-\abs{ab})^2 - D\bigl(s^2 \abs{ab} (\alpha+2)- \abs{ab}(\alpha-2)\bigr).
        \end{aligned}
    \end{equation}
    Therefore,
    \begin{equation}
        \label{eq:h-normalized-calculation}
        \begin{aligned}
            \Bigl ( \frac {h}{1- \omega } \Bigr )^2
             & = \frac {D^2\abs {q_ab+q_ba}^2 }{ \abs {ab}^2} \frac { \abs {ab}}{M^2}                                         \\
             & = \frac{\abs{ab}^2}{M^2}
            \Bigl((\alpha+2)(1-\abs{ab})^2-D\bigl((\alpha+2)s^2-(\alpha-2)\bigr)\Bigr)                                        \\
             & = \frac {( \alpha +2)(1- \abs {ab})^2-D \bigl (( \alpha +2)s^2-( \alpha -2) \bigr )}{( \alpha -2 \abs {ab})^2} \\
             & = \frac {( \alpha +2)(1- \abs {ab})^2-4D-D( \alpha +2)(s^2-1)}{( \alpha -2 \abs {ab})^2}                       \\
             & = \frac {( \alpha -2)(1+ \abs {ab})^2-D( \alpha +2)(s^2-1)}{( \alpha -2 \abs {ab})^2}                          \\
             & = \beta ^2+ \gamma (1-s^2),
        \end{aligned}
    \end{equation}
    where the penultimate equality follows from
    \begin{equation}
        \label{eq:h-constant-factor}
        \begin{aligned}
            ( \alpha +2)(1- \abs {ab})^2-4D
             & =( \alpha +2)(1-2 \abs {ab}+ \abs {ab}^2)
            -4(1+ \abs {ab}^2- \alpha \abs {ab})                              \\
             & =( \alpha -2)+2( \alpha -2) \abs {ab}+( \alpha -2) \abs {ab}^2 \\
             & =( \alpha -2)(1+ \abs {ab})^2.
        \end{aligned}
    \end{equation}

    \par It remains to prove \eqref{eq:h-shape}. The definition of $\gamma$ and \eqref{eq:beta-gamma-bound} give $0\le\gamma\le1$.
    Since $0<\abs{ab}<\varrho\le1$ and $\alpha=\varrho+\varrho^{-1}$, equations \eqref{eq:beta-alpha} and \eqref{eq:beta-gamma-bound} yield
    \begin{equation}
        \varrho-\beta^2-\varrho\gamma
        =\frac{(\alpha-2)(1+\abs{ab})(1+\varrho)(\varrho-\abs{ab})}{(\alpha-2\abs{ab})^2}
        \ge0.
    \end{equation}
    Thus $\beta^2\le(1-\gamma)\varrho$, and \eqref{eq:h-identity} gives the convex-combination bound
    \begin{equation}
        \begin{aligned}
            \left(\frac{h}{1-\omega}\right)^2
             & =\beta^2+\gamma(1-s^2)             \\
             & \le(1-\gamma)\varrho+\gamma(1-s^2)
            \le\max\set{\varrho,1-s^2}.
        \end{aligned}
    \end{equation}
    Taking nonnegative square roots proves \eqref{eq:h-shape}, including the case $\varrho=1$.
\end{proof}

\subsection{The second bound}

\begin{proposition}[Second bound]
    \label{prop:joint2}
    Suppose $s^2 \ge \eta $ and put $ \kappa = \frac {1}{2}(1+ \frac{\eta}{s^2}) \in [ \frac {1}{2},1]$. Then
    \begin{equation}
        \label{eq:C0}
        \Lc \le \mathcal{C}_0:=2 \kappa \, m( \abs {a})+2(1- \kappa ) \, m( \abs {b})+ \sqrt {1- \frac { \eta ^2}{s^2}} \, \sqrt {1-s^2} \; \frac { \J _ \eta (s)}{ \pi }.
    \end{equation}
    Moreover, for fixed $a$ and $b$, the function
    $c\mapsto\mathcal{C}_0(ca,cb)-\Hc(ca,cb)$
    is strictly decreasing for $1\le c<\frac{1}{\abs{a}}$.
\end{proposition}

\begin{proof}
    Direct computation shows that
    \begin{equation}
        \kappa  = \frac { \abs {a-b}^2+ \abs {a}^2- \abs {b}^2}{2 \abs {a-b}^2}
        = \frac {\abs{a}^2 - \re(\bar{a}b)}{ \abs {a-b}^2}
        = \re \frac {a \bar{a} - \bar{a}b}{ \abs {a-b}^2}
        = \re \frac {a}{a-b}.
    \end{equation}
    Put $\tau = \im \frac{a}{a-b}$. Then $ \frac {a}{a-b}= \kappa + \mathrm {i} \tau $. In addition,
    \begin{equation}
        \tau = \im \frac{a}{a-b} = \im \frac{a(\bar{a} - \bar{b})}{ \abs{a-b}^2} = - \frac{\im ( a \bar{b})}{\abs{a-b}^2}.
    \end{equation}
    Since
    \begin{equation}
        \begin{aligned}
            4 \bigl ( \abs {ab}^2-( \re (a \bar {b}))^2 \bigr ) & = \left(2\abs{ab} - 2\re(\bar{a} b)\right)\left(2\abs{ab} + 2\re(\bar{a} b)\right)                          \\
                                                                & = \bigl ( \abs {a-b}^2-( \abs {a}- \abs {b})^2 \bigr ) \bigl (( \abs {a}+ \abs {b})^2- \abs {a-b}^2 \bigr ) \\
                                                                & = (\abs{a}+ \abs{b})^4 (s^2 - \eta^2)(1-s^2),
        \end{aligned}
    \end{equation}
    we have
    \begin{equation} \label{eq:tau-eta-s}
        \abs { \tau }  = \frac { \abs { \im (a \bar {b})}}{ \abs {a-b}^2} = \frac{\sqrt{\abs{a\bar{b}}^2 - \abs{\re (a\bar{b})}^2}}{\abs{a-b}^2}= \frac { \sqrt {(s^2- \eta ^2)(1-s^2)}}{2s^2}.
    \end{equation}
    Since $-\frac{b}{a-b}=1-\frac{a}{a-b}=1-\kappa-\mathrm {i}\tau$, the definition of $W$ in \eqref{eq:functionals} gives
    \begin{equation}
        \label{eq:second-velocity}
        \begin{aligned}
            \frac {v(a)-v(b)}{a-b}
             & =\frac {a}{a-b}\frac {v(a)}{a}-\frac {b}{a-b}\frac {v(b)}{b}        \\
             & =(\kappa+\mathrm {i}\tau)\frac {v(a)}{a}
            +(1-\kappa-\mathrm {i}\tau)\frac {v(b)}{b}                             \\
             & =\kappa\frac {v(a)}{a}+(1-\kappa)\frac {v(b)}{b}+\mathrm {i}\tau W.
        \end{aligned}
    \end{equation}
    Taking real parts and using $\re(\mathrm {i}\tau W)=-\tau\im W$ together with \eqref{eq:functionals}, we get
    \begin{equation}
        \begin{aligned}
            \Lc & = \kappa \re \frac {v(a)}{a}+(1-\kappa) \re \frac {v(b)}{b}+ \re (\mathrm {i}\tau W) \\
                & = \kappa u_a+(1- \kappa )u_b- \tau \im W.
        \end{aligned}
    \end{equation}
    Since $0\le\eta\le s^2$, we have $\frac{1}{2}\le\kappa\le1$, so both weights $\kappa$ and $1-\kappa$ are nonnegative.
    Starting from the preceding identity, we use (R) and (W) from Proposition~\ref{prop:four},
    together with $\abs{\im W}\le\abs{W}$, to obtain
    \begin{equation}
        \label{eq:second-bound-estimate}
        \begin{aligned}
            \Lc & \le \kappa\abs{u_a}+(1-\kappa)\abs{u_b}+\abs{\tau}\abs{\im W} \\
                & \le 2\kappa m(\abs{a})+2(1-\kappa)m(\abs{b})
            +\abs{\tau}\frac{\abs{a-b}\J(a,b)}{\pi}                             \\
                & =2\kappa m(\abs{a})+2(1-\kappa)m(\abs{b})
            +\sqrt{1-\frac{\eta^2}{s^2}}\sqrt{1-s^2}\frac{\J_{\eta}(s)}{\pi}  =\mathcal{C}_0.
        \end{aligned}
    \end{equation}
    For the penultimate equality, substitute \eqref{eq:tau-eta-s},
    $\abs{a-b}=(\abs{a}+\abs{b})s$, and \eqref{eq:jabjeta}:
    \begin{equation}
        \label{eq:second-angular-normalization}
        \begin{aligned}
            \abs{\tau}\frac{\abs{a-b}\J(a,b)}{\pi}
             & =\frac{\sqrt{(s^2-\eta^2)(1-s^2)}}{2s^2}
            \frac{(\abs{a}+\abs{b})s}{\pi}
            \frac{2\J_{\eta}(s)}{\abs{a}+\abs{b}}                                \\
             & =\sqrt{1-\frac{\eta^2}{s^2}}\sqrt{1-s^2}\frac{\J_{\eta}(s)}{\pi}.
        \end{aligned}
    \end{equation}
    \par Finally, under $(a,b)\mapsto(ca,cb)$, the quantities $\varrho,\eta,s,\kappa$ and $\J_{\eta}(s)$ are unchanged, while
    \begin{equation}
        \label{eq:second-target-decomposition}
        \Hc(ca,cb)=\Hc(a,b)-2\log c.
    \end{equation}
    Take $1\le x<y<\frac{1}{\abs{a}}$.
    The angular terms cancel, so
    \begin{equation}
        \label{eq:second-bound-scaling}
        \begin{aligned}
                & \bigl(\mathcal{C}_0(ya,yb)-\Hc(ya,yb)\bigr)
            -\bigl(\mathcal{C}_0(xa,xb)-\Hc(xa,xb)\bigr)                                                                                         \\
            ={} & 2\kappa\bigl(m(y\abs{a})-m(x\abs{a})+\log(\frac{y}{x})\bigr)+2(1-\kappa)\bigl(m(y\abs{b})-m(x\abs{b})+\log(\frac{y}{x})\bigr).
        \end{aligned}
    \end{equation}
    All four radii belong to $(0,1)$. Each parenthesised expression is an increment of $m(q)+\log q$, and is strictly negative by Lemma~\ref{lem:elliptic}(b).
    Since $\kappa\ge\frac{1}{2}>0$ and $1-\kappa\ge0$, their weighted sum is strictly negative, including when $\kappa=1$.
    Hence $c\mapsto\mathcal{C}_0(ca,cb)-\Hc(ca,cb)$ is strictly decreasing on $[1,\frac{1}{\abs{a}})$.
\end{proof}

\section{Proof of the infinitesimal bound}\label{sec:regions}

We first prove $ \Lc \le \Hc $ for $a,b \in \D $ with $a \ne b$ and $ \abs {a} \ge \abs {b}>0$,
and every $ \mu $ with $ \norm { \mu }_ \infty \le 1$.
The nonzero pairs are split into four regions; Table~\ref{tab:regions} summarises the split.

\begin{table}[ht]
    \centering\small
    \caption{The four regions and the estimates used to prove $\Lc\le\Hc$.}\label{tab:regions}
    \renewcommand{\arraystretch}{1.2}
    \begin{tabular}{@{}>{\raggedright\arraybackslash}p{1.05cm}>{\raggedright\arraybackslash}p{5.5cm}>{\raggedright\arraybackslash}p{3cm}>{\raggedright\arraybackslash}p{4.7cm}@{}}
        \toprule
        Region & Conditions                                                                                                                        & Bound for~$\Lc$                     & Main estimates used \\
        \midrule
        I      & $s^2\ge\frac{8}{9}$                                                                                                               & (A) and (R)
               & Lemmas~\ref{lem:elliptic}(b)--(d), \ref{lem:weighted-log}, \ref{lem:large-chord-log}                                                                                                          \\
        II     & $s^2<\frac{8}{9}$, $\beta\le\frac{2}{3}$                                                                                          & Prop.~\ref{prop:joint1}
               & \eqref{eq:J0-bounds}; Lemmas~\ref{lem:coefficient}, \ref{lem:baseline}, \ref{lem:angular-majorant}                                                                                            \\
        III    & $s^2<\frac{8}{9}$, $\beta>\frac{2}{3}$, $s^2\le1-\varrho^4$                                                                       & Prop.~\ref{prop:joint1}
               & \eqref{eq:J0-coarse}; Lemmas~\ref{lem:elliptic}(c), \ref{lem:coefficient}, \ref{lem:shape}                                                                                                    \\
        IV     & $s^2<\frac{8}{9}$, $\beta>\frac{2}{3}$, $s^2>1-\varrho^4$                                                                         & Prop.~\ref{prop:joint2} and~scaling
               & \eqref{eq:J0-bounds}; Lemmas~\ref{lem:elliptic}(c),(d), \ref{lem:simplex}, \ref{lem:region4-g-increasing}, \ref{lem:one-variable}                                                             \\
        \bottomrule
    \end{tabular}
\end{table}

\subsection{Region I}

\begin{proposition}
    \label{prop:region1}
    If $\frac{8}{9} \le s^2 \le 1$, then $ \Lc \le \Hc $.
\end{proposition}

\begin{proof}
    By (A), (R), Lemma~\ref{lem:elliptic}(b), and Lemma~\ref{lem:weighted-log} with $t_1 = \abs {a}$ and $t_2 = \abs {b}$,
    \begin{equation}
        \begin{aligned}
            \Lc  \le \frac { \abs {a}u_a+ \abs {b}u_b}{ \abs {a}+ \abs {b}}+ \A (s) & \le \frac{2\abs{a} m(\abs{a}) + 2 \abs{b} m(\abs{b})}{\abs{a} + \abs{b}} + \A(s)                                   \\
                                                                                    & \le \frac {2 \abs {a} \log (\frac{4}{\abs {a}})+2 \abs {b} \log (\frac{4}{\abs {b}})}{ \abs {a}+ \abs {b}}+ \A (s) \\
                                                                                    & \le 2 \log \frac {8}{ \abs {a}+ \abs {b}}+ \A (s).
        \end{aligned}
    \end{equation}
    Put $c= \sqrt {1-s^2} \in [0, \frac {1}{3}]$. Since
    $p_s \sqrt {1-p_s^2}=\frac{s}{2}$, Lemma~\ref{lem:elliptic}(c),(d) give
    \begin{equation}
        \A (s)=c \, B(p_s) \le 2c \log \frac {8c_0}{s} \le 2c \log \frac {9}{s}.
    \end{equation}
    Since $s= \sqrt {1-c^2}$, Lemma~\ref{lem:large-chord-log} gives
    \begin{equation}
        \A (s) \le 2c \log \frac {9}{s} \le 2 \log \frac {2}{s}.
    \end{equation}
    Combining the two bounds yields $ \Lc \le 2 \log \frac {16}{( \abs {a}+ \abs {b})s}= \Hc $.
\end{proof}

\subsection{Region II}

\begin{lemma}
    \label{lem:baseline}
    With the notation of \eqref{eq:pair} and \eqref{eq:joint-notation}, we have
    \begin{equation}
        \omega B( \delta )+(1- \omega )(2 \log \frac {8}{s}-s^2)< \Hc .
    \end{equation}
\end{lemma}

\begin{proof}
    Put $r=\abs{a}$, $q=\abs{b}$, $S=(r+q)^2$ and $P=\abs{1-\bar{a}b}^2=D+Ss^2$.
    Then $D,M,S>0$, $P>D$ and $0<\omega<1$. The definition of $\delta$ gives $\delta=\frac{\sqrt{S}s}{\sqrt{P}}$ and $1-\delta^2=\frac{D}{P}$, so $0<\delta<1$ and
    \begin{equation}
        \label{eq:baseline-delta-product}
        \delta\sqrt{1-\delta^2}=\frac{\sqrt{S}s\sqrt{D}}{P}.
    \end{equation}
    Lemma~\ref{lem:elliptic}(c) therefore gives
    \begin{equation}
        \label{eq:bdelta}
        \frac{B(\delta)}{2}
        \le\log\frac{4c_0}{\delta\sqrt{1-\delta^2}}
        =\log\frac{4c_0P}{\sqrt{S}s\sqrt{D}}.
    \end{equation}
    Using $\frac{\Hc}{2}=\log(\frac{8}{s})+\log(\frac{2}{\sqrt{S}})$ and $s^2=\frac{P-D}{S}$, we obtain
    \begin{equation}
        \label{eq:baseline-reduction}
        \begin{aligned}
                  & \frac{1}{2}\left(\omega B(\delta)+(1-\omega)\left(2\log\frac{8}{s}-s^2\right)-\Hc\right)           \\
            ={}   & \frac{\omega B(\delta)}{2}-\omega\log\frac{8}{s}+\log\frac{\sqrt{S}}{2}-\frac{(1-\omega)(P-D)}{2S} \\
            \le{} & (1-\omega)\log\frac{\sqrt{S}}{2}+\omega\log\frac{c_0P}{4\sqrt{D}}-\frac{(1-\omega)(P-D)}{2S}.
        \end{aligned}
    \end{equation}
    The arithmetic--geometric mean inequality and $2(r^2+q^2)\ge S$ give
    \begin{equation}
        \label{eq:baseline-radial-bound}
        \begin{aligned}
            M
             & =\frac{r^2+q^2}{2}\bigl((1-r^2)+(1-q^2)\bigr)+\frac{(r^2-q^2)^2}{2} \\
             & \ge(r^2+q^2)\sqrt{D}
            \ge\frac{S\sqrt{D}}{2}.
        \end{aligned}
    \end{equation}
    Applying $\log z\le z-1$ to $z=\frac{MP}{2SD}>0$ yields
    \begin{equation}
        \label{eq:baseline-log-tangent}
        \log\frac{c_0P}{4\sqrt{D}}
        \le\log\frac{c_0S\sqrt{D}}{2M}+\frac{MP}{2SD}-1.
    \end{equation}
    Substitute this into \eqref{eq:baseline-reduction} and use $\omega M=(1-\omega)D$ from \eqref{eq:Z}.
    Since $\sqrt{S}=r+q<2$ and $S-M=2rq(1+rq)>0$, the radial bound \eqref{eq:baseline-radial-bound} and Lemma~\ref{lem:elliptic}(d) give
    \begin{equation}
        \label{eq:baseline-maximum}
        \begin{aligned}
                  & \frac{1}{2}\left(\omega B(\delta)+(1-\omega)\left(2\log\frac{8}{s}-s^2\right)-\Hc\right) \\
            \le{} & (1-\omega)\log\frac{\sqrt{S}}{2}
            +\omega\left(\log\frac{c_0S\sqrt{D}}{2M}-1+\frac{M}{2S}\right)                                   \\
            <{}   & \omega\left(\log c_0-\frac{1}{2}\right)<0.
        \end{aligned}
    \end{equation}
    This proves the lemma.
\end{proof}

\begingroup
\renewcommand{\thetheorem}{\getrefnumber{lem:angular-majorant}}
\begin{lemma}
    \label{main:angular}
    For $0<t \le \frac{2 \sqrt {2}}{3}$,
    \begin{equation}
        \tag{\ref*{eq:angular-majorant}}
        \label{eq:angular-majorant-restated}
        \sqrt {1-\frac{5t^2}{9}} \, \bigl (3+(2 \log \frac {8}{t}-3)(1-t^2)^{-\frac{1}{8}} \bigr ) \le 2 \log \frac {8}{t}-t^2.
    \end{equation}
\end{lemma}
\addtocounter{theorem}{-1}
\endgroup

The proof is elementary calculus and is given in \hyperref[app:angular]{Appendix~\ref*{app:elementary}}.

\begin{proposition}
    \label{prop:region2}
    If $0<s^2<\frac{8}{9}$ and $ \beta \le \frac{2}{3}$, then $ \Lc < \Hc $.
\end{proposition}

\begin{proof}
    Using \eqref{eq:h-identity}, \eqref{eq:beta-gamma-bound}, and $0\le\beta\le\frac{2}{3}$, we obtain
    \begin{equation}
        \Bigl ( \frac {h}{1- \omega } \Bigr )^2=( \beta ^2+ \gamma )(1-s^2)+ \beta ^2s^2 \le (1-s^2)+ \frac {4}{9}s^2=1- \frac {5s^2}{9},
    \end{equation}
    which yields
    \begin{equation} \label{eq:h-bound}
        h\le(1-\omega)\sqrt {1-\frac{5s^2}{9}}.
    \end{equation}
    Then we have
    \begin{equation}
        \label{eq:region2-L-bound}
        \begin{aligned}
            \Lc & \le\omega B(\delta)+h\frac {\Jz(s)}{\pi}                                     \\
                & \le\omega B(\delta)+(1-\omega)\sqrt {1-\frac {5s^2}{9}}\,\frac {\Jz(s)}{\pi} \\
                & \le\omega B(\delta)+(1-\omega)\sqrt {1-\frac {5s^2}{9}}\,
            \bigl(3+(2\log\frac {8}{s}-3)(1-s^2)^{-\frac{1}{8}}\bigr)                          \\
                & \le\omega B(\delta)+(1-\omega)\bigl(2\log\frac {8}{s}-s^2\bigr).
        \end{aligned}
    \end{equation}
    The first inequality follows from Proposition~\ref{prop:joint1}, and the second uses \eqref{eq:h-bound}.
    The third follows from \eqref{eq:J0-bounds}, and the last inequality follows from Lemma~\ref{lem:angular-majorant}.
    Finally, Lemma~\ref{lem:baseline} shows that the last expression in \eqref{eq:region2-L-bound} is strictly less than $\Hc$.
    This proves $\Lc<\Hc$.
\end{proof}

\subsection{Region III}

\begingroup
\renewcommand{\thetheorem}{\getrefnumber{lem:shape}}
\begin{lemma}
    \label{main:shape}
    If $0< \abs {b} \le \abs {a}<1$ and $0 \le s^2 \le 1- \varrho ^4$, then
    \begin{equation}
        \tag{\ref*{eq:shape-estimate}}
        \label{eq:shape-estimate-restated}
        \omega \log \frac {c_0 \bigl (D+( \abs {a}+ \abs {b})^2s^2 \bigr )}{2( \abs {a}+ \abs {b}) \sqrt {D}}< \log \frac {2}{ \abs {a}+ \abs {b}}.
    \end{equation}
\end{lemma}
\addtocounter{theorem}{-1}
\endgroup

The proof is a case analysis given in \hyperref[app:shape]{Appendix~\ref*{app:elementary}}.

\begin{proposition}
    \label{prop:region3}
    If $0<s^2<1$ and $s^2 \le 1- \varrho ^4$, then $ \Lc < \Hc $.
\end{proposition}

\begin{proof}
    Under the hypothesis, both $ \sqrt {1-s^2}$ and $ \sqrt { \varrho }$ are at most $(1-s^2)^{\frac{1}{8}}$. So \eqref{eq:h-shape}
    and \eqref{eq:J0-coarse} give $\frac{h \Jz (s)}{\pi} \le (1- \omega ) \, 2 \log \frac {8}{s}$, and
    Proposition~\ref{prop:joint1} gives $ \Lc \le 2 \log \frac {8}{s}+ \omega \bigl (B( \delta )-2 \log \frac {8}{s} \bigr )$. By
    Lemma~\ref{lem:elliptic}(c) and $ \delta \sqrt {1- \delta ^2}= \frac{\abs {a-b} \sqrt {D}}{D+ \abs {a-b}^2}$,
    \begin{equation}
        B( \delta )-2 \log \frac {8}{s} \le 2 \log \frac {c_0 \bigl (D+( \abs {a}+ \abs {b})^2s^2 \bigr )}{2( \abs {a}+ \abs {b}) \sqrt {D}}.
    \end{equation}
    Lemma~\ref{lem:shape} gives $ \Lc <2 \log \frac {8}{s}+2 \log \frac {2}{ \abs {a}+ \abs {b}}= \Hc $.
\end{proof}

\subsection{Region IV}

\begingroup
\renewcommand{\thetheorem}{\getrefnumber{lem:one-variable}}
\begin{lemma}
    \label{main:one-variable}
    For $ \frac {1}{ \sqrt {3}} \le t <1$, we have
    \begin{equation}
        \tag{\ref*{eq:one-variable}}
        \label{eq:one-variable-restated}
        \frac {3(1- t )}{2} \log \frac {9}{ \sqrt {\frac{3(1- t )}{2}}}
        \le (1- t ^{\frac{3}{2}}) \, \log \frac {8}{ \sqrt {1- t ^4}}+ \frac {3}{2} \bigl ( t ^{\frac{3}{2}}- t ^2 \bigr )+ \frac {3}{2} \log \frac {2}{1+ t }.
    \end{equation}
\end{lemma}
\addtocounter{theorem}{-1}
\endgroup

The proof is given in \hyperref[app:one-variable]{Appendix~\ref*{app:elementary}}.

\begin{proposition}
    \label{prop:region4}
    If $ \beta >\frac{2}{3}$ and $1- \varrho ^4<s^2<\frac{8}{9}$, then $\mathcal{C}_0< \Hc $, and hence $ \Lc < \Hc $.
\end{proposition}

\begin{proof}
    In this proof, $f$ denotes an auxiliary real function defined below.
    The global quasiconformal map $f$ is not used here, so this local reuse of notation causes no ambiguity.

    From $1-\varrho^4<\frac{8}{9}$ we get $\varrho>\frac{1}{\sqrt{3}}$,
    while $\beta>0$ forces $\varrho<1$. Since
    $0<\eta<1-\varrho<1-\varrho^4<s^2$, Proposition~\ref{prop:joint2}
    applies, with $\frac{1}{2}<\kappa<1$.

    Put $x=\abs{a}^2$ and $x_0=\frac{1+\varrho}{2}$. By \eqref{eq:beta},
    \begin{equation}
        \label{eq:beta-v}
        \beta=\frac{\sqrt{\varrho}(1-\varrho)(1+x\varrho)}
        {1+(1-2x)\varrho^2}.
    \end{equation}
    If $x\le x_0$, the numerator is at most its value at $x_0$,
    while the denominator is at least $1-\varrho^3>0$. Hence
    \begin{equation}
        \begin{aligned}
            \beta \le\frac{\sqrt{\varrho}(\varrho^2+\varrho+2)}{2(1+\varrho+\varrho^2)}
             & <\frac{(1+\varrho)(\varrho^2+\varrho+2)} {4(1+\varrho+\varrho^2)} \\
             & =\frac{2}{3}
            -\frac{(1-\varrho)(3\varrho^2+\varrho+2)}
             {12(1+\varrho+\varrho^2)}
            <\frac{2}{3},
        \end{aligned}
    \end{equation}
    where the first strict inequality uses
    $\sqrt{\varrho}<\frac{1+\varrho}{2}$.
    This contradicts $\beta>\frac{2}{3}$, so $x>x_0$.
    By the scaling monotonicity in Proposition~\ref{prop:joint2}, it suffices
    to prove the desired inequality at the threshold radii
    \begin{equation}
        \label{eq:threshold-radii}
        r_a=\sqrt{x_0},\qquad r_b=\varrho r_a,
    \end{equation}
    with the same $\varrho$ and $s$.
    In the remainder of the proof, $\mathcal{C}_0$ and $\Hc$ are evaluated
    at this threshold scale.

    Since $\varrho^2>\frac{1}{3}$ and $\varrho>\frac{1}{2}$, we have
    \begin{equation}
        r_b^2=\frac{\varrho^2(1+\varrho)}{2}
        >\frac{(\frac{1}{3})(\frac{3}{2})}{2}=\frac{1}{4}.
    \end{equation}
    Thus $\frac{1}{2}<r_b<r_a<1$.
    Define $f(q)=q\log(\frac{9}{\sqrt{q}})$ for $0<q\le1$.
    For $u\in\set{r_a,r_b}$, Lemma~\ref{lem:elliptic}(c),(d) gives
    \begin{equation}
        \label{eq:threshold-m-majorant}
        \begin{aligned}
            m(u)
             & =\frac{1-u^2}{2}B(u)                     \\
             & \le(1-u^2)\log\frac{4c_0}{u\sqrt{1-u^2}}
            \le f(1-u^2),
        \end{aligned}
    \end{equation}
    because $\frac{4c_0}{u}<8c_0<9$.
    Moreover, $f'(q)=\log(\frac{9}{\sqrt{q}})-\frac{1}{2}>0$ and
    $f''(q)=-\frac{1}{2q}<0$, where the positivity of $f'$ follows
    from $\log9>\frac{1}{2}$.
    Hence $f$ is increasing and strictly concave on $(0,1]$.

    Put
    \begin{equation}
        A_1=\kappa m(r_a)+(1-\kappa)m(r_b).
    \end{equation}
    Since $f(1-r_a^2)\le f(1-r_b^2)$ and $\kappa\ge\frac{1}{2}$,
    removing the excess weight on the smaller term and then applying
    Jensen's inequality gives
    \begin{equation}
        \label{eq:threshold-radial}
        \begin{aligned}
            A_1
             & \le\kappa f(1-r_a^2)+(1-\kappa)f(1-r_b^2) \\
             & \le\frac{f(1-r_a^2)+f(1-r_b^2)}{2}        \\
             & \le f\left(1-\frac{r_a^2+r_b^2}{2}\right)
            \le f\left(\frac{3(1-\varrho)}{2}\right).
        \end{aligned}
    \end{equation}
    The last inequality uses the monotonicity of $f$ and
    \begin{equation}
        1-\frac{r_a^2+r_b^2}{2}
        =\frac{1-\varrho}{4}(\varrho^2+2\varrho+3)
        \le\frac{3(1-\varrho)}{2}<\frac{3}{4}<1.
    \end{equation}
    The final bound in \eqref{eq:threshold-radial} depends only on
    $\varrho$, and is therefore uniform in $s$ and $\kappa$.

    Since $0<\eta<s^2<1$, Lemma~\ref{lem:simplex} yields
    \begin{equation}
        \sqrt{1-\frac{\eta^2}{s^2}}\,\J_{\eta}(s)
        \le\sqrt{1-\eta^2}\,\J_{\eta}(s)\le\Jz(s).
    \end{equation}
    Let $g$ be the function in Lemma~\ref{lem:region4-g-increasing}.
    Substituting the preceding estimate and \eqref{eq:J0-bounds}
    into \eqref{eq:C0}, and using the definition of $g$, gives
    \begin{equation}
        \label{eq:threshold-angular}
        \begin{aligned}
            \frac{\mathcal{C}_0}{2}
             & \le A_1+\frac{\sqrt{1-s^2}}{2}\frac{\Jz(s)}{\pi}            \\
             & \le A_1+\frac{3}{2}(1-s^2)^{\frac{1}{2}}
            +\left(\log\frac{8}{s}-\frac{3}{2}\right)(1-s^2)^{\frac{3}{8}} \\
             & =A_1+\log\frac{8}{s}-g(s^2).
        \end{aligned}
    \end{equation}
    At the threshold scale, $r_a+r_b=(1+\varrho)\sqrt{\frac{1+\varrho}{2}}$,
    so
    \begin{equation}
        \frac{\Hc}{2}
        =\log\frac{8}{s}+\log\frac{2}{r_a+r_b}
        =\log\frac{8}{s}+\frac{3}{2}\log\frac{2}{1+\varrho}.
    \end{equation}
    Because $\frac{1}{\sqrt{3}}<\varrho<1$, Lemma~\ref{lem:one-variable}
    applies with $t=\varrho$. Together with \eqref{eq:threshold-radial},
    it gives
    \begin{equation}
        A_1\le f\left(\frac{3(1-\varrho)}{2}\right)
        \le g(1-\varrho^4)+\frac{3}{2}\log\frac{2}{1+\varrho}.
    \end{equation}
    Subtracting the expression for $\frac{\Hc}{2}$ from
    \eqref{eq:threshold-angular}, we obtain
    \begin{equation}
        \mathcal{C}_0-\Hc
        \le2\bigl(g(1-\varrho^4)-g(s^2)\bigr)<0.
    \end{equation}
    The last inequality follows from $0<1-\varrho^4<s^2<\frac{8}{9}$
    and the strict monotonicity in Lemma~\ref{lem:region4-g-increasing}.
    The original pair is a common enlargement of the threshold pair,
    so Proposition~\ref{prop:joint2} gives $\mathcal{C}_0<\Hc$
    at the original pair as well.
\end{proof}

\subsection{Conclusion}

\begin{proof}[Proof of Theorem~\ref{thm:kernel}]
    Fix $ \mu \in L^ \infty ( \D )$ with $ \norm { \mu }_ \infty \le 1$ and distinct $a,b \in \D $.
    By symmetry, label the points so that $ \abs {a} \ge \abs {b}$; then $a \ne 0$.
    First suppose $b \ne 0$. The four alternatives
    \begin{enumerate}[label=(\Roman*)]
        \item $s^2 \ge \frac{8}{9}$;
        \item $s^2<\frac{8}{9}$ and $ \beta \le \frac{2}{3}$;
        \item $s^2<\frac{8}{9}$, $ \beta >\frac{2}{3}$ and $s^2 \le 1- \varrho ^4$;
        \item $s^2<\frac{8}{9}$, $ \beta >\frac{2}{3}$ and $s^2>1- \varrho ^4$
    \end{enumerate}
    cover all cases. Propositions~\ref{prop:region1}, \ref{prop:region2}, \ref{prop:region3} and \ref{prop:region4}
    give $ \Lc _ \mu (a,b) \le \Hc (a,b)$, proving \eqref{eq:kernel} for nonzero pairs.

    It remains to consider $b=0$. Since $a \ne 0$ and $v(0)=0$, we have
    $ \Lc _ \mu (a,0)= \re (\frac{v(a)}{a})$.
    The radial conclusion of Theorem~\ref{thm:schwarz} gives
    $ \abs {f_t(a)} \le \varphi _{K_t}( \abs {a})$ for $0 \le t<1$, with equality at $t=0$.
    Differentiating the logarithms at $t=0^+$ and using $K_t'(0)=2$ and Lemma~\ref{lem:elliptic}(a),(b), we obtain
    \begin{equation}
        \Lc _ \mu (a,0) \le 2m( \abs {a}) \le 2 \log \frac {4}{ \abs {a}}
        <2 \log \frac {16}{ \abs {a}}= \Hc (a,0).
    \end{equation}
\end{proof}

\appendix
\section{Elementary identities and inequalities}\label{app:elementary}\label{app:elementary-identities}

\begin{lemma}
    \label{lem:log-modulus-derivative}
    Let $I\subset\R$ be an interval and let $f:I\to\C$ be differentiable at an interior point $t_0\in I$.
    If $f(t_0)\ne0$, then
    \begin{equation}
        \label{eq:log-modulus-derivative}
        \left.\frac{\dd}{\dd t}\log\abs{f(t)}\right|_{t=t_0}
        =\frac{\re\left(f'(t_0)\overline{f(t_0)}\right)}{\abs{f(t_0)}^2}
        =\re\frac{f'(t_0)}{f(t_0)}.
    \end{equation}
    The same formula holds for $t_0 \in \partial I$ if $f$ has the corresponding one-sided derivative there and $f(t_0)\ne0$.
\end{lemma}

\begin{proof}
    By continuity, $f$ is nonzero near $t_0$ within $I$.
    Since complex conjugation commutes with differentiation in the real variable $t$, the product rule gives
    \begin{equation}
        \label{eq:modulus-squared-derivative}
        \left.\frac{\dd}{\dd t}\abs{f(t)}^2\right|_{t=t_0}
        =\left.\frac{\dd}{\dd t}\left(f(t)\overline{f(t)}\right)\right|_{t=t_0}
        =2\re\left(f'(t_0)\overline{f(t_0)}\right).
    \end{equation}
    Apply the chain rule to $\log\abs{f(t)}=\frac{1}{2}\log\abs{f(t)}^2$ and use
    $\frac{\overline{f(t_0)}}{\abs{f(t_0)}^2}=\frac{1}{f(t_0)}$ to obtain \eqref{eq:log-modulus-derivative}.
    The same argument applies to one-sided derivatives.
\end{proof}

\begin{lemma}
    \label{lem:log-series-bounds}
    Let $N \ge 1$ be an integer and $0<y<1$. Then
    \begin{equation}
        \label{eq:log-series-bounds}
        2 \sum _{n=0}^{N-1} \frac {y^{2n+1}}{2n+1}
        < \log \frac {1+y}{1-y}
        <2 \sum _{n=0}^{N-1} \frac {y^{2n+1}}{2n+1}
        + \frac {2y^{2N+1}}{(2N+1)(1-y^2)}.
    \end{equation}
\end{lemma}

\begin{proof}
    The first inequality follows directly from \eqref{eq:log-series}. For the second inequality, replacing the denominator of each item in the tail by $2N+1$ yields
    \begin{equation}
        \label{eq:log-series-tail}
        0<2 \sum _{n=N}^{ \infty } \frac {y^{2n+1}}{2n+1}
        < \frac {2}{2N+1} \sum _{n=N}^{ \infty }y^{2n+1}
        = \frac {2y^{2N+1}}{(2N+1)(1-y^2)}.
    \end{equation}
    Adding the first $N$ terms proves the stated bounds.
\end{proof}

All numerical estimates for logarithms used in this paper are collected in the following lemma.

\begin{lemma}
    \label{lem:rational-logs}
    The following numerical bounds hold:
    \begin{equation}
        \label{eq:numerical-log-bounds}
        \begin{aligned}
             & \frac{7}{60}<\frac{2}{17}<\log\frac{9}{8},\quad
            \log\frac{3}{2}<\frac{73}{180}<\frac{5}{12},                         \\[2pt]
             & \frac{69}{100}<\frac{842}{1215}<\log2<\frac{25}{36}<\frac{7}{10}, \\[2pt]
             & \frac{13}{12}<\log3<\frac{11}{10},\quad
            \log\frac{27}{8}<\frac{73}{60}<\frac{5}{4},                          \\[2pt]
             & \log8>2>\frac{3}{2}>\frac{1}{2},\quad
            \log9\ge\frac{8}{9}>\frac{1}{2},\quad
            \log72>\frac{21}{5}.
        \end{aligned}
    \end{equation}
\end{lemma}

\begin{proof}
    Apply Lemma~\ref{lem:log-series-bounds} with $y=\frac{1}{17}$ and $N=1$.
    The lower bound gives
    \begin{equation}
        \frac {7}{60}< \frac {2}{17}< \log \frac {9}{8}.
    \end{equation}
    With $y=\frac{1}{5}$ and $N=1$, the upper bound gives
    \begin{equation}
        \log\frac{3}{2}<\frac{73}{180}<\frac{5}{12}.
    \end{equation}
    With $y=\frac{1}{3}$, the lower bound at $N=2$ and the upper bound at $N=1$ give
    \begin{equation}
        \frac{69}{100}<2\left(\frac{1}{3}+\frac{1}{81}\right)<\log2
        <\frac{25}{36}<\frac{7}{10}.
    \end{equation}
    At $y=\frac{1}{3}$ and $N=3$, the lower bound gives the sharper estimate
    \begin{equation}
        \label{eq:log2-sharper-lower}
        \frac{69}{100}<\frac{842}{1215}<\log2.
    \end{equation}
    With $y=\frac{1}{2}$ and $N=2$, Lemma~\ref{lem:log-series-bounds} gives
    \begin{equation}
        \frac{13}{12}<\log3<\frac{11}{10}.
    \end{equation}
    The bound for $\log(\frac{3}{2})$ implies
    \begin{equation}
        \label{eq:one-variable-log-twenty-seven-eighths}
        \log\frac{27}{8}=3\log\frac{3}{2}
        <\frac{73}{60}<\frac{5}{4}.
    \end{equation}
    The lower bound for $\log2$ yields
    \begin{equation}
        \label{eq:one-variable-log8}
        \log8=3\log2>3\cdot\frac{69}{100}>2>\frac{3}{2}>\frac{1}{2}.
    \end{equation}
    Also,
    \begin{equation}
        \label{eq:log9-lower}
        \log9=\int_{1}^{9}\frac{\dd u}{u}\ge\frac{8}{9}>\frac{1}{2}>\frac{1}{4}.
    \end{equation}
    Finally, the lower bounds for $\log2$ and $\log3$ give
    \begin{equation}
        \label{eq:log72-lower}
        \log72=2\log3+3\log2
        >\frac{13}{6}+\frac{207}{100}>\frac{21}{5}.
    \end{equation}
\end{proof}

The Feynman-parameter identity represents a product of reciprocal powers as a parameter integral with a single combined denominator.
For an introduction to Feynman parametrization, see \cite{KCL}.
Lemma~\ref{lem:feynman} follows directly from \cite[Eq.~(53)]{KCL}; we include a proof for the reader's convenience.

\begin{lemma}[Feynman-parameter identity]
    \label{lem:feynman}
    Let $A,B,C>0$, and set
    \begin{equation}
        \Sigma = \set {( \xi , \zeta ) \in \R ^2: \xi >0, \  \zeta >0, \  \xi + \zeta <1}, \qquad \ell =1- \xi - \zeta .
    \end{equation}
    Then
    \begin{equation}
        \label{eq:feynman}
        \frac {1}{ \sqrt {ABC}}= \frac {1}{2 \pi } \int _ \Sigma \frac { \dd \xi \dd \zeta }{ \sqrt { \xi \zeta \ell } \, ( \xi A+ \zeta B+ \ell C)^{\frac{3}{2}}}.
    \end{equation}
\end{lemma}

\begin{proof}
    For $ \alpha , \lambda >0$, the substitution $u= \lambda t$ in the gamma integral \cite[8.310.1]{GR} gives
    \begin{equation}
        \int _0^ \infty t^{ \alpha -1} \ee ^{- \lambda t} \dd t
        = \lambda ^{- \alpha } \int _0^ \infty u^{ \alpha -1} \ee ^{-u} \dd u
        = \frac { \Gamma ( \alpha )}{ \lambda ^ \alpha }.
    \end{equation}
    Apply this with $ \alpha = \frac {1}{2}$ and $ \lambda =A,B,C$, and multiply the three identities.
    Since the integrands are nonnegative, Tonelli's theorem gives
    \begin{equation}
        \frac {1}{ \sqrt {ABC}}
        = \frac {1}{ \Gamma ( \frac {1}{2})^3}
        \int _{(0, \infty )^3} \frac { \ee ^{-Au_1-Bu_2-Cu_3}}{ \sqrt {u_1u_2u_3}} \dd u_1 \dd u_2 \dd u_3.
    \end{equation}
    Separate the total scale from the relative weights by setting
    \begin{equation}
        u_1=r \xi , \qquad u_2=r \zeta , \qquad u_3=r \ell ,
        \qquad \ell =1- \xi - \zeta , \quad r>0, \quad ( \xi , \zeta ) \in \Sigma .
    \end{equation}
    This is a smooth bijection from $(0, \infty ) \times \Sigma $ onto $(0, \infty )^3$, with smooth inverse
    $r=u_1+u_2+u_3$, $ \xi =\frac{u_1}{r}$, $ \zeta =\frac{u_2}{r}$. Its Jacobian determinant is
    \begin{equation}
        \det
        \begin{pmatrix}
            \xi   & r  & 0  \\
            \zeta & 0  & r  \\
            \ell  & -r & -r
        \end{pmatrix}
        =r^2( \xi + \zeta + \ell )=r^2.
    \end{equation}
    Thus $ \dd u_1 \dd u_2 \dd u_3=r^2 \dd r \dd \xi \dd \zeta $, whereas
    $ \sqrt {u_1u_2u_3}=r^{\frac{3}{2}} \sqrt { \xi \zeta \ell }$.
    After this change of variables, Tonelli's theorem allows us to integrate first in $r$.
    Since $ \xi A+ \zeta B+ \ell C \ge \min \set {A,B,C}>0$, the gamma integral with
    $ \alpha = \frac {3}{2}$ yields
    \begin{equation}
        \begin{aligned}
            \frac {1}{ \sqrt {ABC}}
             & = \frac {1}{ \Gamma ( \frac {1}{2})^3}
            \int _ \Sigma \frac {1}{ \sqrt { \xi \zeta \ell }}
            \left ( \int _0^ \infty r^{\frac{1}{2}} \ee ^{-r( \xi A+ \zeta B+ \ell C)} \dd r \right ) \dd \xi \dd \zeta \\
             & = \frac { \Gamma ( \frac {3}{2})}{ \Gamma ( \frac {1}{2})^3}
            \int _ \Sigma \frac { \dd \xi \dd \zeta }{ \sqrt { \xi \zeta \ell } \, ( \xi A+ \zeta B+ \ell C)^{\frac{3}{2}}}.
        \end{aligned}
    \end{equation}
    All integrals are nonnegative, and the initial product of gamma integrals is finite; hence this argument
    also proves convergence of the simplex integral despite its singularities at the boundary.
    Finally, $ \Gamma ( \frac {1}{2})= \sqrt { \pi }$ and
    $ \Gamma ( \frac {3}{2})= \frac {1}{2} \Gamma ( \frac {1}{2})$ \cite[8.338.2, 8.331.1]{GR}, so
    \begin{equation}
        \frac { \Gamma ( \frac {3}{2})}{ \Gamma ( \frac {1}{2})^3}
        = \frac { \frac{\sqrt{\pi}}{2}}{( \sqrt { \pi })^3}= \frac {1}{2 \pi },
    \end{equation}
    which proves \eqref{eq:feynman}.
\end{proof}

\begin{lemma}[A planar radial integral]
    \label{lem:plane-radial-integral}
    Let $c \in \C $ and $ \Delta >0$. Then
    \begin{equation}
        \label{eq:plane-radial-integral}
        \int _ \C \frac { \dA (z)}{( \abs {z-c}^2+ \Delta )^{\frac{3}{2}}}= \frac {2 \pi }{ \sqrt { \Delta }}.
    \end{equation}
\end{lemma}

\begin{proof}
    The integrand is nonnegative, so the polar change of variables
    $z=c+r \ee ^{ \mathrm {i} \theta }$, with area element $r \dd r \dd \theta $, gives
    \begin{equation}
        \begin{aligned}
            \int _ \C \frac { \dA (z)}{( \abs {z-c}^2+ \Delta )^{\frac{3}{2}}}
             & = \int _0^{2 \pi } \int _0^ \infty \frac {r}{(r^2+ \Delta )^{\frac{3}{2}}} \dd r \dd \theta                              \\
             & =2 \pi \int _0^ \infty \frac {r}{(r^2+ \Delta )^{\frac{3}{2}}} \dd r = \pi \int _ \Delta ^ \infty u^{-\frac{3}{2}} \dd u \\
             & = \frac {2 \pi }{ \sqrt { \Delta }}.
        \end{aligned}
    \end{equation}
\end{proof}

\begin{lemma}
    \label{lem:weighted-log}
    Let $t_1,t_2 \in (0,1]$. Then
    \begin{equation}
        \label{eq:weighted-log}
        \frac {t_1 \log (\frac{4}{t_1})+t_2 \log (\frac{4}{t_2})}{t_1+t_2} \le \log \frac {8}{t_1+t_2}.
    \end{equation}
    Equality holds if and only if $t_1=t_2$.
\end{lemma}
\begin{proof}
    The function $t\mapsto\log t$ is strictly concave on $(0,\infty)$,
    since its second derivative is $-\frac{1}{t^2}<0$.
    The weights $\frac{t_1}{t_1+t_2}$ and $\frac{t_2}{t_1+t_2}$ are positive and sum to one.
    Applying Jensen's inequality to the positive numbers $\frac{4}{t_1}$ and $\frac{4}{t_2}$, we obtain
    \begin{equation}
        \label{eq:weighted-log-jensen}
        \begin{aligned}
            \frac{t_1\log(\frac{4}{t_1})+t_2\log(\frac{4}{t_2})}{t_1+t_2}
             & =\frac{t_1}{t_1+t_2}\log\frac{4}{t_1}
            +\frac{t_2}{t_1+t_2}\log\frac{4}{t_2}            \\
             & \le\log\left(\frac{t_1}{t_1+t_2}\frac{4}{t_1}
            +\frac{t_2}{t_1+t_2}\frac{4}{t_2}\right)         \\
             & =\log\frac{8}{t_1+t_2}.
        \end{aligned}
    \end{equation}
    This proves \eqref{eq:weighted-log}. Since the logarithm is strictly concave and both weights are positive,
    equality holds if and only if $\frac{4}{t_1}=\frac{4}{t_2}$, or equivalently, $t_1=t_2$.
\end{proof}

\begin{lemma}
    \label{lem:large-chord-log}
    Let $t \in [0,\frac{1}{3}]$. Then
    \begin{equation}
        \label{eq:large-chord-log}
        t \log \frac {9}{ \sqrt {1-t^2}} \le \log \frac {2}{ \sqrt {1-t^2}}.
    \end{equation}
    Equality holds if and only if $t=\frac{1}{3}$.
\end{lemma}

\begin{proof}
    For $0 \le t \le \frac{1}{3}$, put
    \begin{equation}
        f(t)= \log \frac {2}{ \sqrt {1-t^2}}-t \log \frac {9}{ \sqrt {1-t^2}}.
    \end{equation}
    For $0<t<\frac{1}{3}$, differentiation gives
    \begin{equation}
        f'(t)= \frac {t}{1+t}- \log \frac {9}{ \sqrt {1-t^2}} \le \frac {1}{4}- \log 9<0.
    \end{equation}
    Thus $f$ is strictly decreasing on $[0,\frac{1}{3}]$.
    At the right endpoint,
    \begin{equation}
        f \Bigl ( \frac {1}{3} \Bigr )= \log \frac {3}{ \sqrt {2}}- \frac {1}{3} \log \Bigl [ \Bigl ( \frac {3}{ \sqrt {2}} \Bigr )^3 \Bigr ]=0.
    \end{equation}
    Therefore $f(t) \ge f(\frac{1}{3})=0$, with equality if and only if $t=\frac{1}{3}$.
\end{proof}

\begin{lemma}
    \label{lem:angular-majorant}
    \label{app:angular}
    For $0<t \le \frac{2\sqrt{2}}{3}$,
    \begin{equation}
        \label{eq:angular-majorant}
        \sqrt {1-\frac{5t^2}{9}} \, \bigl (3+(2 \log \frac {8}{t}-3)(1-t^2)^{-\frac{1}{8}} \bigr ) \le 2 \log \frac {8}{t}-t^2.
    \end{equation}
\end{lemma}
\begin{proof}
    For $0<t\le\frac{2\sqrt{2}}{3}$ and $y\in\mathbb{R}$, set
    \begin{equation}
        \label{eq:angular-two-variable}
        f(t,y)=\sqrt{1-\frac{5t^2}{9}}
        \left(3+(y-3)(1-t^2)^{-\frac{1}{8}}\right)-y+t^2.
    \end{equation}
    Writing $1-t^2=(1-\frac{9t^2}{8})\cdot 1 + \frac{9t^2}{8} \cdot \frac{1}{9}$, and applying Jensen's inequality to the concave function $u\mapsto u^{\frac{1}{4}}$ at the endpoints $u=\frac{1}{9}$ and $u=1$, with weight $1-\frac{9t^2}{8}$ assigned to $u=1$,
    we obtain
    \begin{equation}
        \begin{aligned}
            (1-t^2)^{\frac{1}{4}}
             & \ge\left(1-\frac{9t^2}{8}\right)\cdot1^{\frac{1}{4}}
            +\frac{9t^2}{8}\left(\frac{1}{9}\right)^{\frac{1}{4}}   \\
             & =1-\frac{9t^2}{8}+\frac{9t^2}{8\sqrt{3}}
            >1-\frac{5t^2}{9},
        \end{aligned}
    \end{equation}
    where the last inequality is strict because $t^2>0$ and $\frac{1}{\sqrt{3}}>\frac{41}{81}$.
    Consequently,
    \begin{equation}
        \label{eq:angular-y-derivative}
        \frac{\partial f}{\partial y}(t,y)
        =\sqrt{1-\frac{5t^2}{9}}(1-t^2)^{-\frac{1}{8}}-1<0.
    \end{equation}
    Since
    \begin{equation}
        \label{eq:angular-log-lower}
        2\log\frac{8}{t}\ge\log72>\frac{21}{5},
    \end{equation}
    we have
    \begin{equation}
        \label{eq:angular-y-reduction}
        f\left(t,2\log\frac{8}{t}\right)
        \le f\left(t,\frac{21}{5}\right).
    \end{equation}
    Thus it suffices to prove
    \begin{equation}
        \label{eq:fixed}
        \sqrt{1-\frac{5t^2}{9}}
        \left(3+\frac{6}{5}(1-t^2)^{-\frac{1}{8}}\right)
        \le\frac{21}{5}-t^2.
    \end{equation}
    \par Put $x=(1-t^2)^{\frac{1}{8}}$, $x \in [(\frac{1}{3})^{\frac{1}{4}},1] \subset [\frac{3}{4}, 1]$.
    Under this substitution, \eqref{eq:fixed} is equivalent to
    $\varphi(x)\ge0$, where
    \begin{equation}
        \varphi(x)=\log\frac{x(16+5x^8)}
        {(5x+2)\sqrt{4+5x^8}}
        \qquad (\frac{3}{4}\le x\le1).
    \end{equation}
    Differentiation yields
    \begin{equation}
        \label{eq:angular-derivative}
        \varphi'(x)
        =-\frac{2x^7g(x)}{(5x+2)(4+5x^8)(16+5x^8)},
    \end{equation}
    where
    \begin{equation}
        g(x)=60+400x-125x^8-250x^9-64x^{-8}.
    \end{equation}
    The function $g$ is strictly concave on $(0,\infty)$, because
    \begin{equation}
        g''(x)=-7000x^6-18000x^7-4608x^{-10}<0.
    \end{equation}
    At the endpoints,
    \begin{equation}
        g\left(\frac{3}{4}\right)
        =360-\frac{625}{2}\left(\frac{3}{4}\right)^8
        -64\left(\frac{4}{3}\right)^8<360-576<0,
    \end{equation}
    and $g(1)=21>0$.
    Hence $g$ has exactly one zero in $(\frac{3}{4},1)$, changing sign from negative to positive.
    It follows from \eqref{eq:angular-derivative} that $\varphi$ first increases
    and then decreases, so its minimum on $[\frac{3}{4},1]$ occurs at an endpoint.
    We have $\varphi(1)=0$. At $x=\frac{3}{4}$, we have $x^8<\frac{1}{9}$ and
    \begin{equation}
        (16+5x^8)^2-\left(5+\frac{2}{x}\right)^2(4+5x^8)
        =\frac{225x^{16}-1205x^8+188}{9}
        >\frac{487}{81}>0.
    \end{equation}
    Since the quantities being squared are positive, this proves
    $\varphi(\frac{3}{4})>0$. Therefore $\varphi(x)\ge0$ on $[\frac{3}{4},1]$,
    which proves \eqref{eq:fixed}. Combining this with
    \eqref{eq:angular-y-reduction} proves the lemma.
\end{proof}

\begin{lemma}
    \label{lem:shape}
    \label{app:shape}
    If $0< \abs {b} \le \abs {a}<1$ and $0 \le s^2 \le 1- \varrho ^4$, then
    \begin{equation}
        \label{eq:shape-estimate}
        \omega \log \frac {c_0 \bigl (D+( \abs {a}+ \abs {b})^2s^2 \bigr )}{2( \abs {a}+ \abs {b}) \sqrt {D}}< \log \frac {2}{ \abs {a}+ \abs {b}}.
    \end{equation}
\end{lemma}

\begin{proof}
    Put
    \begin{equation}
        \xi=\frac{1+\abs{ab}}{1-\abs{ab}}>1,\qquad
        \zeta=1-\left(\frac{\abs{a}-\abs{b}}{1-\abs{ab}}\right)^2\in(0,1].
    \end{equation}
    Direct computation gives
    \begin{equation}
        \begin{aligned}
            D          & =(1-\abs{ab})^2\zeta,\qquad
                       & \sqrt{\xi^2-\zeta}          & =\frac{\abs{a}+\abs{b}}{1-\abs{ab}}>0, \\
            1-\abs{ab} & =\frac{2}{\xi+1},\qquad
                       & \omega                      & =\frac{\zeta}{\xi}.
        \end{aligned}
    \end{equation}
    Hence, with
    \begin{equation}
        T=\frac{c_0\bigl(D+(\abs{a}+\abs{b})^2s^2\bigr)}
        {2(\abs{a}+\abs{b})\sqrt{D}}>0,
    \end{equation}
    the desired conclusion is equivalent to
    \begin{equation}
        \label{eq:shape-target}
        \zeta\log T<\xi\log\frac{\xi+1}{\sqrt{\xi^2-\zeta}}.
    \end{equation}
    By assumption and $\varrho=\frac{1-\eta}{1+\eta}$,
    \begin{equation}
        s^2\le1-\varrho^4
        =\frac{8\eta(1+\eta^2)}{(1+\eta)^4}
        \le8\eta
        =\frac{8(1-\abs{ab})\sqrt{1-\zeta}}{\abs{a}+\abs{b}}.
    \end{equation}
    Substituting this into the definition of $T$, we obtain
    \begin{equation}
        \label{eq:shape-T-upper}
        \frac{T}{c_0}
        =\frac{\sqrt{D}}{2(\abs{a}+\abs{b})}
        +\frac{(\abs{a}+\abs{b})s^2}{2\sqrt{D}}
        \le\frac{\sqrt{\zeta}}{2\sqrt{\xi^2-\zeta}}+4\sqrt{\frac{1-\zeta}{\zeta}}.
    \end{equation}

    \begin{enumerate}[label=(\roman*),leftmargin=2em]
        \item Suppose $0<\zeta\le\frac{8}{9}$.
              Since $\xi^2-\zeta\ge\xi^2(1-\zeta)>0$, $\frac{\sqrt{1-\zeta}}{\sqrt{\xi^2-\zeta}} \le \frac{1}{\xi}$.
              Together with $\frac{\zeta}{1-\zeta}\le8$, we have $\frac{\zeta}{8(\xi^2-\zeta)} < \frac{1}{\xi^2}$.
              Therefore, \eqref{eq:shape-T-upper} yields
              \begin{equation}
                  \begin{aligned}
                      \frac{T^2}{16c_0^2}
                       & \le\left(\sqrt{\frac{1-\zeta}{\zeta}}+\frac{\sqrt{\zeta}}{8\sqrt{\xi^2-\zeta}}\right)^2          \\
                       & =\frac{1-\zeta}{\zeta} +\frac{\sqrt{1-\zeta}}{4\sqrt{\xi^2-\zeta}}+\frac{\zeta}{64(\xi^2-\zeta)} \\
                       & \le\frac{1-\zeta}{\zeta}+\frac{1}{4\xi}+\frac{1}{8\xi^2}
                      \le\frac{1-\zeta}{\zeta}+\frac{3}{8\xi}.
                  \end{aligned}
              \end{equation}
              Applying $\log u\le u-1$ to
              $u=\frac{4}{3}\left(\frac{1-\zeta}{\zeta}+\frac{3}{8\xi}\right)>0$ gives
              \begin{equation}
                  \label{eq:shape-log-upper}
                  \begin{aligned}
                      \log T^2 & \le\log\left(12c_0^2 \cdot \frac{4}{3}\left(\frac{1-\zeta}{\zeta}+\frac{3}{8\xi}\right)\right) \\
                               & =\log(12c_0^2) +\log\left(\frac{4}{3}\left(\frac{1-\zeta}{\zeta}+\frac{3}{8\xi}\right)\right)  \\
                               & \le\log(12c_0^2)+\frac{4}{3\zeta}-\frac{7}{3}+\frac{1}{2\xi}.
                  \end{aligned}
              \end{equation}

              For $y\ge0$, the function $\log(1+y)-\frac{2y}{2+y}$ vanishes at $0$
              and has derivative $\frac{y^2}{(1+y)(2+y)^2}\ge0$.
              Thus $\log(1+y)\ge\frac{2y}{2+y}$.
              Together with $-\log(1-u)\ge u$ for $0\le u<1$, this implies
              \begin{equation}
                  \label{eq:shape-log-lower}
                  \begin{aligned}
                      \xi\log\frac{\xi+1}{\sqrt{\xi^2-\zeta}}
                       & =\xi\log\left(1+\frac{1}{\xi}\right)
                      -\frac{\xi}{2}\log\left(1-\frac{\zeta}{\xi^2}\right) \\
                       & \ge\frac{2\xi}{2\xi+1}+\frac{\zeta}{2\xi}.
                  \end{aligned}
              \end{equation}
              Since $\frac{4}{3}-\frac{4\xi}{2\xi+1}<0$ and $\frac{1}{\zeta}\ge\frac{9}{8}$,
              multiplying \eqref{eq:shape-log-lower} by $\frac{2}{\zeta}>0$ and subtracting it
              from \eqref{eq:shape-log-upper} gives
              \begin{equation}
                  \begin{aligned}
                      \frac{2}{\zeta}\left(\zeta\log T-\xi\log\frac{\xi+1}{\sqrt{\xi^2-\zeta}}\right)
                       & \le\log(12c_0^2)+\frac{4}{3\zeta}-\frac{7}{3}+\frac{1}{2\xi} - \frac{2}{\zeta}\left(\frac{2\xi}{2\xi+1}+\frac{\zeta}{2\xi}\right) \\
                       & = \log(12c_0^2)-\frac{7}{3}-\frac{1}{2\xi} +\frac{1}{\zeta}\left(\frac{4}{3}-\frac{4\xi}{2\xi+1}\right)                           \\
                       & \le\log(12c_0^2)-\frac{7}{3}-\frac{1}{2\xi}
                      +\frac{9}{8}\left(\frac{4}{3}-\frac{4\xi}{2\xi+1}\right)                                                                             \\
                       & =\log(12c_0^2)-\frac{17}{6}
                      -\frac{(\xi-1)^2}{2\xi(2\xi+1)}.
                  \end{aligned}
              \end{equation}
              Lemma~\ref{lem:c0}\textup{(d)} gives
              \begin{equation}
                  \log(12c_0^2)=2\log c_0+2\log2+\log3
                  <\frac{7}{30}+\frac{25}{18}+\frac{11}{10}
                  <\frac{17}{6}.
              \end{equation}
              Therefore the preceding difference is strictly negative, proving
              \eqref{eq:shape-target} in this case.

        \item Suppose $\frac{8}{9}<\zeta\le1$. Then $\frac{1-\zeta}{\zeta}<\frac{1}{8}$.
              Multiplying \eqref{eq:shape-T-upper} by $\sqrt{\xi^2-\zeta}$
              and applying the Cauchy--Schwarz inequality, we obtain
              \begin{equation}
                  \begin{aligned}
                      T\sqrt{\xi^2-\zeta} & \le c_0\left(\frac{1}{2}\cdot \sqrt{\zeta}
                      +4\sqrt{\frac{1-\zeta}{\zeta}} \cdot \sqrt{\xi^2-\zeta}\right)                                                            \\
                                          & \le c_0\sqrt{\left(\frac{1}{4}+16\frac{1-\zeta}{\zeta}\right)\left( \zeta + (\xi^2 - \zeta)\right)}
                      \le\frac{3c_0}{2}\xi<2\xi.
                  \end{aligned}
              \end{equation}
              Therefore,
              \begin{equation}
                  \begin{aligned}
                      T<\frac{2\xi}{\sqrt{\xi^2-\zeta}}
                       & <\frac{2\xi}{\sqrt{\xi^2-\zeta}}+\frac{(\xi-1)(\xi-\sqrt{\xi^2-\zeta})}{\sqrt{\xi^2-\zeta}} \\
                       & =1+\xi\left(\frac{\xi+1}{\sqrt{\xi^2-\zeta}}-1\right)
                      \le\left(\frac{\xi+1}{\sqrt{\xi^2-\zeta}}\right)^\xi,
                  \end{aligned}
              \end{equation}
              where the last inequality uses the elementary inequality $1+ \xi(u-1) \le u^{\xi}$ for $u \ge 0$ and $\xi \ge 1$.
              Since $0<\sqrt{\xi^2-\zeta}<\xi$, $\frac{\xi+1}{\sqrt{\xi^2-\zeta}} > 1$.
              Taking logarithms and using $0<\zeta\le1$ and $\log(\frac{\xi+1}{\sqrt{\xi^2-\zeta}})>0$, we conclude that
              \begin{equation}
                  \zeta\log T<\zeta\xi\log\frac{\xi+1}{\sqrt{\xi^2-\zeta}}
                  \le\xi\log\frac{\xi+1}{\sqrt{\xi^2-\zeta}},
              \end{equation}
              which proves \eqref{eq:shape-target}.
    \end{enumerate}
\end{proof}

\begin{lemma}
    \label{lem:region4-g-increasing}
    The function $g:(0,\frac{8}{9}]\to\R$ defined by
    \begin{equation}
        \label{eq:region4-G}
        g(t)=\log\frac{8}{\sqrt{t}}\bigl(1-(1-t)^{\frac{3}{8}}\bigr)
        +\frac{3}{2}\bigl((1-t)^{\frac{3}{8}}-(1-t)^{\frac{1}{2}}\bigr)
    \end{equation}
    is strictly increasing.
\end{lemma}

\begin{proof}
    Set $h(t)=1-(1-t)^{\frac{3}{8}}$ for $0\le t<1$.
    Its derivatives are $h'(t)=\frac{3}{8}(1-t)^{-\frac{5}{8}}>0$ and
    $h''(t)=\frac{15}{64}(1-t)^{-\frac{13}{8}}>0$.
    Since $h'$ is increasing and $h(0)=0$, we have
    \begin{equation}
        h(t)=\int_0^t h'(u)\dd u
        \le\int_0^t h'(t)\dd u=th'(t).
    \end{equation}
    Rewriting \eqref{eq:region4-G} as
    \begin{equation}
        g(t)=\left(\log\frac{8}{\sqrt{t}}-\frac{3}{2}\right)h(t)
        +\frac{3}{2}\left(1-\sqrt{1-t}\right),
    \end{equation}
    and differentiating, we obtain
    \begin{equation}
        \begin{aligned}
            g'(t)
             & =-\frac{h(t)}{2t}
            +\left(\log\frac{8}{\sqrt{t}}-\frac{3}{2}\right)h'(t)
            +\frac{3}{4\sqrt{1-t}}                           \\
             & \ge\left(\log\frac{8}{\sqrt{t}}-2\right)h'(t)
            +\frac{3}{4\sqrt{1-t}}>0
            \qquad(0<t<\frac{8}{9}),
        \end{aligned}
    \end{equation}
    where the first inequality uses $\frac{h(t)}{t} \le h'(t)$, and the second uses $h'(t)>0$ and $\log(\frac{8}{\sqrt{t}})>\log8>2$.
\end{proof}

\begin{lemma}
    \label{lem:one-variable}
    \label{app:one-variable}
    For $ \frac {1}{ \sqrt {3}} \le t <1$, we have
    \begin{equation}
        \label{eq:one-variable}
        \frac {3(1- t )}{2} \log \frac {9}{ \sqrt {3(1-t)/2}} \le  (1- t ^{\frac{3}{2}}) \, \log \frac {8}{ \sqrt {1- t ^4}} + \frac {3}{2} \bigl ( t ^{\frac{3}{2}}- t ^2 \bigr )+ \frac {3}{2} \log \frac {2}{1+ t }.
    \end{equation}
\end{lemma}
\begin{proof}
    The conclusion is equivalent to
    \begin{equation}
        \label{eq:one-variable2}
        \frac {3}{2} \log \frac {9}{ \sqrt {3(1-t)/2}} -  \frac{1- t ^{\frac{3}{2}}}{1- t} \, \log \frac {8}{ \sqrt {1- t ^4}} \le \frac{3}{2(1-t)} \left( t ^{\frac{3}{2}}- t ^2 +  \log \frac {2}{1+ t }\right).
    \end{equation}
    We first estimate the left-hand side of \eqref{eq:one-variable2}. Since $t^{\frac{3}{2}}\le\frac{t+t^2}{2}$, $\frac{1-t^{\frac{3}{2}}}{1-t}\ge1+\frac{t}{2}$.
    \begin{equation}
        \begin{aligned}
            \text{LHS of}\,\, \eqref{eq:one-variable2} \le & \frac{3}{2}\log\frac{9}{\sqrt{3(1-t)/2}} - \left( \frac{3}{2} - \frac{1-t}{2}\right)\log\frac{8}{\sqrt{1-t^4}} \\
            =                                              & \frac{3}{2}\log \sqrt{\frac{27P}{32}} + \frac{1-t}{2}\log\frac{8}{\sqrt{1-t^4}},
        \end{aligned}
    \end{equation}
    where $P=(1+t)(1+t^2)$. Direct computation shows that $2<P < 4$ and $4-P=(1-t)(t^2+2t+3)>\frac{13}{3}(1-t)$.
    Applying $\log u\le u-1$ to $u=\frac{P}{4}$ and using
    $\log\frac{27}{8}<\frac{5}{4}$, we obtain
    \begin{equation}
        \label{eq:one-variable-logP-uniform}
        \begin{aligned}
            \frac{3}{4}\log\frac{27P}{32} =\frac{3}{4}\log\frac{27}{8}+\frac{3}{4}\log\frac{P}{4} & < \frac{15}{16}+ \frac{3}{4}\left(\frac{P}{4}-1\right) \\
                                                                                                  & <\frac{15}{16}-\frac{13}{16}(1-t).
        \end{aligned}
    \end{equation}
    Applying the same inequality to $u=\frac{1}{2(1-t)}$ and using
    $\log2<\frac{7}{10}$ gives
    \begin{equation}
        \label{eq:one-variable-weighted-log-uniform}
        \begin{aligned}
            \frac{1-t}{2}\log\frac{8}{\sqrt{2(1-t)}} = & \frac{3}{2}(1-t)\log2
            +\frac{1-t}{4}\log\frac{1}{2(1-t)}                                                                                          \\
            \le                                        & \frac{3}{2}(1-t)\log2+\frac{1}{8}-\frac{1-t}{4} <\frac{1}{8}+\frac{4}{5}(1-t).
        \end{aligned}
    \end{equation}
    Therefore,
    \begin{equation}
        \label{eq:one-variable-uniform-upper}
        \begin{aligned}
            \text{LHS of}\,\, \eqref{eq:one-variable2} \le & \frac{3}{4}\log\frac{27P}{32}
            +\frac{1-t}{2}\log\frac{8}{\sqrt{1-t^4}}                                                     \\
            <                                              & \frac{3}{4}\log\frac{27P}{32}
            +\frac{1-t}{2}\log\frac{8}{\sqrt{2(1-t)}}                                                    \\
            <                                              & \frac{17}{16}-\frac{1-t}{80}<\frac{17}{16}.
        \end{aligned}
    \end{equation}
    On the other hand,
    $\log\frac{2}{1+t}=-\log\left(1-\frac{1-t}{2}\right)\ge\frac{1-t}{2}$.
    Also, $t\ge\frac{1}{\sqrt{3}}>\frac{9}{16}$ gives $\frac{3}{4}<\sqrt{t}<1$.
    Therefore,
    \begin{equation}
        \label{eq:one-variable-remainder-lower}
        \begin{aligned}
            \text{RHS of}\,\, \eqref{eq:one-variable2} \ge & \frac{3}{2(1-t)} (t^{\frac{3}{2}}- t ^2 +  \frac{1-t}{2})                                                               \\
            =                                              & \frac{3t^{\frac{3}{2}}}{2(1+\sqrt{t})}+\frac{3}{4} >  \frac{3}{4}\left(\frac{3}{4}\right)^3+\frac{3}{4} >\frac{17}{16}.
        \end{aligned}
    \end{equation}
    Thus, \eqref{eq:one-variable2} and consequently \eqref{eq:one-variable} are proved.
\end{proof}

\clearpage
\section{Table of notation}\label{app:notation}

\begingroup
\small
\renewcommand{\arraystretch}{1.04}
\begin{longtable}{@{}>{$}p{3.4cm}<{$}>{\raggedright\arraybackslash}p{9.15cm}>{\raggedright\arraybackslash}p{2.2cm}@{}}
    \caption{Table of notation.}\label{tab:notation}                                                                                                                                                                                                       \\
    \toprule
    \text {Symbol}                             & Meaning                                                                                                                                  & Defined in                                                     \\
    \midrule
    \endfirsthead
    \multicolumn{3}{c}{\tablename\ \thetable: Table of notation (continued).}                                                                                                                                                                              \\
    \toprule
    \text {Symbol}                             & Meaning                                                                                                                                  & Defined in                                                     \\
    \midrule
    \endhead
    \midrule
    \multicolumn{3}{r}{Continued on the next page.}                                                                                                                                                                                                        \\
    \endfoot
    \bottomrule
    \endlastfoot
    \mathcal {F}_K                             & orientation-preserving $K$-quasiconformal homeomorphisms of $ \D $ onto itself fixing $0$                                                & \S\ref{subsec:result}                                          \\
    \mu _f                                     & Beltrami coefficient $\frac{f_{\bar{z}}}{f_z}$ of $f$                                                                                    & \S\ref{subsec:qc}                                              \\
    M(K)                                       & Mori constant                                                                                                                            & \eqref{eq:M}                                                   \\
    k, \mu                                     & $\frac{K-1}{K+1}$ and $\frac{\mu_f}{k}$: path endpoint and normalised coefficient                                                        & \eqref{eq:normalized-beltrami-direction}                       \\
    K_t                                        & $\frac{1+t}{1-t}$, the quasiconformal dilatation bound along the path                                                                    & \S\ref{subsec:strategy}                                        \\
    f_z,\ f_{\bar{z}},\ J_f                    & weak complex derivatives and Jacobian                                                                                                    & \S\ref{subsec:qc}                                              \\
    f^{\mu}                                    & normalised solution with Beltrami coefficient $\mu$                                                                                      & \S\ref{subsec:beltrami}                                        \\
    \nu,\ \dot{\lambda}                        & variation directions in the original and image planes                                                                                    & Thm.~\ref{thm:variation}, Cor.~\ref{cor:transported-variation} \\
    F_t                                        & $f^{t\mu^{\mathrm{e}}}$, with $f_t=F_t|_{\D}$                                                                                            & \S\ref{sec:reduction}                                          \\
    V[ \nu ]                                   & normalised first-variation operator on the sphere                                                                                        & \eqref{eq:V}                                                   \\
    q'= \sqrt {1-q^2}                          & complementary modulus                                                                                                                    & \S\ref{subsec:elliptic}                                        \\
    \K , \  \mathcal{G} , \  \varphi _K        & complete elliptic integral, Grötzsch modulus, distortion function                                                                        & \eqref{eq:elliptic}                                            \\
    B(q),\ m(q),\ c_0                          & $\frac{4}{\pi}\K(q)\K(q')$, $\frac{q'^2}{2}B(q)$, and $\frac{1}{8}\exp\bigl(\frac{2}{\pi}\K(\frac{1}{\sqrt{2}})^2\bigr)\approx1.1152$    & \eqref{eq:Bm}                                                  \\
    \delta(a,b), \  \delta                     & pseudo-hyperbolic distance $\frac{\abs{a-b}}{\abs{1-\bar{a}b}}$                                                                          & \S\ref{subsec:distortion}                                      \\
    d_h(a,b)                                   & hyperbolic distance on $ \D $, with density $\frac{2}{1-\abs{z}^2}$                                                                      & \S\ref{subsec:distortion}                                      \\
    \Theta                                     & angular distortion parameter; for $(0,a,b)$, $ \sin ( \frac{\Theta}{2})=s$                                                               & Thm.~\ref{thm:angle}                                           \\
    ( \alpha )_n                               & rising factorial $ \alpha ( \alpha +1) \cdots ( \alpha +n-1)$, with $( \alpha )_0=1$                                                     & \S\ref{subsec:series}                                          \\
    b_n, \  d_n                                & $ \binom {2n}{n}4^{-n}$ and $H_{2n}-H_n$ ($H_n$ harmonic numbers)                                                                        & \eqref{eq:bn}                                                  \\
    \iota                                      & reflection in the unit circle, $ \iota (z)=\frac{1}{\bar{z}}$                                                                            & \S\ref{subsec:reflection}                                      \\
    \mu , \  \mu ^{ \mathrm {e}}, \  f_t, \  v & deformation direction, its reflection, deformed disk maps, velocity                                                                      & \S\ref{sec:reduction}                                          \\
    \Lc , \  \Lc _ \mu (a,b)                   & chord log-speed $ \re \frac {v(a)-v(b)}{a-b}$; the subscript records the deformation direction                                           & \eqref{eq:L-kernel}, \eqref{eq:functionals}                    \\
    \Hc , \  \Hc (a,b)                         & the target $2 \log (\frac{16}{\abs{a-b}})$                                                                                               & \eqref{eq:L-kernel}                                            \\
    s,\ \eta,\ q_a,\ q_b                       & relative chord, radial imbalance, and radial weights $q_a,q_b$                                                                           & \eqref{eq:pair}                                                \\
    u_a, \  u_b                                & radial log-speeds $ \re (\frac{v(a)}{a})$ and $ \re (\frac{v(b)}{b})$                                                                    & \eqref{eq:functionals}                                         \\
    W                                          & complex quantity $\frac{v(a)}{a}-\frac{v(b)}{b}$                                                                                         & \eqref{eq:functionals}                                         \\
    p_s                                        & angular parameter $ \sin ( \frac {1}{2} \arcsin s)$                                                                                      & \eqref{eq:Aang}                                                \\
    \A (s)                                     & angular majorant $ \sqrt {1-s^2} \, B(p_s)$                                                                                              & \eqref{eq:Aang}                                                \\
    \textup{(R), (S), (A), (W)}                & radial, two-point Schwarz, angular, and three-pole integral bounds                                                                       & \eqref{eq:R}--\eqref{eq:W}                                     \\
    \J , \  \J _ \eta , \  \Jz                 & three-pole integral and its normalisations                                                                                               & \S\ref{sec:threepole}                                          \\
    \Sigma ; \  \xi , \zeta , \ell             & simplex $ \set {( \xi , \zeta ): \xi >0, \  \zeta >0, \  \xi + \zeta <1}$ and its barycentric coordinates, with $ \ell =1- \xi - \zeta $ & Prop.~\ref{prop:simplex-representation}                        \\
    \Theta_0                                   & denominator term $u(1-u)+u^2(1-q^2)s^2$ at $\eta=0$                                                                                      & Lem.~\ref{lem:simplex}                                         \\
    \varrho                                    & $ \frac{\abs{b}}{\abs{a}} \in (0,1]$                                                                                                     & \eqref{eq:joint-notation}                                      \\
    D, \  M, \  \omega                         & $(1- \abs {a}^2)(1- \abs {b}^2)$, $ \abs {a}^2+ \abs {b}^2-2 \abs {ab}^2$, $\frac{D}{1-\abs{ab}^2}$                                      & \eqref{eq:joint-notation}                                      \\
    h, \  \beta                                & coefficient of the first bound and balance parameter                                                                                     & \eqref{eq:h}, \eqref{eq:beta}                                  \\
    \gamma                                     & coefficient of $1-s^2$ in $(\frac{h}{1-\omega})^2= \beta ^2+ \gamma (1-s^2)$                                                             & \eqref{eq:beta}                                                \\
    \kappa                                     & weight $ \frac {1}{2}(1+ \frac{\eta}{s^2})$ in the second bound                                                                          & Prop.~\ref{prop:joint2}                                        \\
    \mathcal{C}_0                              & the second upper bound for $ \Lc $                                                                                                       & \eqref{eq:C0}                                                  \\
\end{longtable}
\endgroup

\clearpage
\printbibliography

\end{document}